\documentclass[11pt,letterpaper]{article}

\usepackage{amsmath,amssymb,enumitem}
\usepackage{titlesec}
\titleformat{\subsection}[runin]
{\normalfont\large\bfseries}{\thesubsection}{0.5em}{}
\titlespacing{\subsection}{0pt}{*0.5}{*1}
\titlespacing{\section}{0pt}{*0.8}{*0.8}
\usepackage{comment,tikz}
\usepackage{graphicx}
\usepackage{setspace}
\usepackage{hyperref}
\hypersetup{colorlinks=true, citecolor=blue, linkcolor=blue, urlcolor=blue}
\usepackage{cleveref}
\usepackage[T1]{fontenc}
\usepackage{float,graphicx,comment,color,lmodern,cite}
\usepackage[margin=1in]{geometry}
\usepackage{amsmath,amssymb,amsthm,mathtools}
\usepackage{microtype,enumitem,booktabs,array,longtable,needspace}
\usepackage{hyperref}
\setlist[enumerate]{label=\textup{(\arabic*)},ref=(\arabic*),leftmargin=*,itemsep=.3em}
\setlist[itemize]{leftmargin=*,itemsep=.25em}
\numberwithin{equation}{section}
\newtheorem{theorem}{Theorem}[section]
\newtheorem{lemma}[theorem]{Lemma}
\newtheorem{proposition}[theorem]{Proposition}
\newtheorem{corollary}[theorem]{Corollary}
\newtheorem{conjecture}{Conjecture}

\theoremstyle{definition}

\theoremstyle{remark}

\newcommand{\PP}{\mathbb P}
\newcommand{\EE}{\mathbb E}

\newcommand{\tr}{\operatorname{tr}}
\newcommand{\head}{\operatorname{head}}
\newcommand{\tail}{\operatorname{tail}}

\tikzset{
  v/.style={circle, fill=black, inner sep=1.7pt},
  e/.style={thick},
  de/.style={thick}
}

\makeatletter
\def\@settitle{\begin{center}%
  \baselineskip14\p@\relax
    \normalfont\LARGE
  \@title
  \vspace{0.5in}
  \end{center}%
}
\makeatother

\makeatletter
\renewcommand\part{\@startsection{part}{0}\z@
  {2\linespacing\@plus\linespacing}{\linespacing}
  {\normalfont\Large\bfseries\centering\def\@secnumfont{\bfseries}}}
\makeatother

\author{}
\date{}

\begin{document}

\title{Counting Eulerian Orientations \\ {\large in} \\ Regular Graphs}

\author{Srinjoy Srimani \and Jacques Verstra\"ete\thanks{Research supported by NSF Award DMS-2347832.} \and Jason Yang}
\date{\small Department of Mathematics, University of California, San Diego\\
      \small La Jolla, CA 92093, United States\\
      \small \texttt{ssrimani@ucsd.edu}, \texttt{jverstraete@math.ucsd.edu}, \texttt{jay049@ucsd.edu}}

\maketitle

\begin{abstract}
Let $G$ be a pseudograph, and let $\varepsilon(G)$ denote the number of Eulerian orientations of $G$. Lieb and Wu, and later Schrijver, showed that if $G$ has degree sequence $(d_1,d_2,\dots,d_n)$, where each $d_i$ is a positive even integer, then  
\[ \varepsilon(G) \geq \prod_{i = 1}^n 2^{-d_i/2} \binom{d_i}{d_i/2}.\]
Schrijver conjectured that for simple $d$-regular graphs, $\varepsilon(G)\leq\varepsilon(K_{d+1})^{n/(d+1)}$, with equality attained by disjoint unions of cliques of order $d+1$. In this paper, we prove a multigraph bound that implies the known case $d=4$, and we also improve Schrijver's lower bound for $d$-regular graphs, using an inductive relation between $\varepsilon(G)$ and the number $\ell_i(G)$ of cycles of each possible length $i \geq 1$ in $G$. We define 
\[ \ell(G) = \sum_{i = 1}^n \frac{\ell_i(G)}{(d - 1)^i}\]
and we prove with $\lambda = (d - 1)\log(1 + 1/(d - 1))$ that for all $d \geq 2$, 
\[ \varepsilon(G) \geq \exp(\lambda \ell(G)) \cdot S_d(n),\]
where $S_d(n)$ denotes $2^{-dn/2}\binom{d}{d/2}^n$. The bound sharply improves the lower bound of Lieb, Wu and Schrijver, and gives good estimates for $\varepsilon(G)$ for various graphs, such as integer lattices. The bound is almost tight for random regular graphs, and we conjecture for each even $d \geq 4$ that there exists a constant $c_d > 0$ such that $\varepsilon(G) \leq \exp(c_d\lambda \ell(G)) \cdot S_d(n)$ for every $n$-vertex $d$-regular pseudograph.
\end{abstract}

\section{Introduction}

A {\em pseudograph} is a pair $(V,E)$ where $E$ is a multiset 
of unordered pairs of elements of $V$. We allow pairs of the form $\{a,a\}$, which we refer to as {\em loops}. In words, both loops and multiple edges are permitted. A {\em multigraph} is a pseudograph without loops. Parallel edge copies are distinct, and a loop has two possible directions, according to which incidence is its tail. An \emph{Eulerian orientation} of a pseudograph is an assignment of arrows to its edges such that the in-degree and out-degree at every vertex are equal. By Euler's Theorem, such an orientation exists if and only if all the degrees are even. Let $\varepsilon(G)$ denote the number of Eulerian orientations of a pseudograph $G$.

The topic of Eulerian orientations and Eulerian tours is one of the fundamental research subjects in graph theory, with applications to algorithms, optimization, probability, statistical physics, and network theory.
The number of Eulerian orientations of a $4$-regular undirected graph is $(-1)^{|V(G)|+c(G)}T_G(0,-2)$, where $T_G$ is the Tutte polynomial and $c(G)$ is the number of components -- see Jaeger, Vertigan and Welsh~\cite{JVW90,Wel93}. Enumeration of orientations is also connected to the study of statistical physics through Pauling's model for the residual entropy of ice~\cite{Pau35} and Lieb's study of square ice~\cite{Lie67}. The well-known BEST theorem -- see Tutte and Smith~\cite{TS41} and van Aardenne-Ehrenfest and de Bruijn~\cite{AB51} -- states that the number of Eulerian tours in an Eulerian digraph $G$ with distinct arcs, connected on its nonisolated vertices, with a fixed first arc leaving $w$, is given by the formula
\[t_w(G) \cdot \prod_{v:\,d^+(v)>0} (d^+(v) - 1)!,\]
where $t_w(G)$ is the number of spanning arborescences on the nonisolated vertices directed towards $w$. The quantity $t_w(G)$ is computed using a matrix-tree theorem for directed graphs, which can be achieved efficiently for a directed graph.  The problem of computing the number of Eulerian tours in undirected graphs is $\#\mathrm{P}$-complete -- see Brightwell and Winkler~\cite{BW05}. Counting Eulerian orientations is also $\#\mathrm{P}$-complete, as shown by Mihail and Winkler~\cite{MW96}.
In this paper, we establish new lower bounds on the number of Eulerian orientations of $d$-regular pseudographs and an upper bound for multigraphs with degrees two and four. 

\subsection{Lower bounds.} 
 Lieb and Wu~\cite{LW72}, and later Schrijver~\cite{Sch83}, proved that for an undirected pseudograph $G$ on $n$ vertices with degree sequence $(d_1,d_2,\dots,d_n)$ where $d_i$ is even, 
\begin{equation}\label{eq:Schrijverbound} \varepsilon(G) \geq \prod_{i=1}^n 2^{-d_i/2} \binom{d_i}{d_i/2}.
\end{equation}
 In the $d$-regular case, (\ref{eq:Schrijverbound}) 
becomes 
\begin{equation}\label{eq:Schrijverboundregular} 
\varepsilon(G) \geq 2^{-dn/2}\binom{d}{d/2}^n \eqqcolon S_d(n).
\end{equation}
Schrijver~\cite{Sch83} showed using a counting argument that for each fixed even $d\geq4$, there exist a constant $c_d>0$ and $d$-regular pseudographs $G$ for which $\varepsilon(G)\leq c_d\sqrt{n}\cdot S_d(n)$. For simple graphs the same bound follows from Creed~\cite[Theorem~3.49]{Creed10}. Our first theorem gives a lower bound for $\varepsilon(G)$ which is sensitive to the number of cycles of various lengths in the pseudograph $G$. We let $\ell_i(G)$ denote the number of simple cycles of length $i$ in a pseudograph $G$, counted up to rotation and reversal. Thus $\ell_1(G)$ is the number of loops, and $t$ parallel edges between two distinct vertices contribute $\binom{t}{2}$ to $\ell_2(G)$. 

\begin{theorem}\label{thm:main}
    Let $n \geq 1$, and let $d \geq 2$ be even. Then for any $n$-vertex $d$-regular pseudograph $G$,   
    \begin{equation} \label{eq:cyclebound}
    \varepsilon(G) \geq  \prod_{i = 1}^n \Bigl(1 + \frac{1}{d - 1}\Bigr)^{\frac{\ell_i(G)}{(d - 1)^{i - 1}}} \cdot S_d(n).
\end{equation}
\end{theorem}

We may write the product in \eqref{eq:cyclebound} as 
$$\exp\Bigl(\lambda \sum_{i = 1}^n \frac{\ell_i(G)}{(d - 1)^i}\Bigr)$$
where $\lambda = (d - 1)\log(1 + 1/(d - 1))$ and $\log$ denotes the natural logarithm. Theorem \ref{thm:main} is applied effectively 
when we have good bounds on the sum in the exponent. The following theorem will quantify this statement. We shall refer to 
\begin{equation}\label{eq:cyclesum}
\ell(G) := \sum_{i = 1}^n \frac{\ell_i(G)}{(d - 1)^i}
\end{equation}
as the {\em weighted cycle sum of $G$.}
We let $w_i(G)$ denote the number of rooted, directed closed walks of length $i$ in $G$ with no immediate reversal. Each walk is a sequence of directed edge copies with a distinguished first edge, and the reversal condition also applies from the last edge to the first. For example, if $G$ is a cycle of length $i$, then $\ell_i(G) = 1$ whereas $w_i(G) = 2i$. We shall refer to 
\begin{equation}\label{eq:walksum}
w(G) := \sum_{i = 1}^n \frac{w_i(G)}{2i(d - 1)^i}
\end{equation}
as the {\em weighted walk sum of $G$.}

\begin{theorem}\label{thm:walks}
There exists an absolute constant $c$ such that for any simple $n$-vertex $d$-regular graph $G$ with $d\ge3$, the following hold. 
\begin{enumerate}[labelwidth=2.5em, leftmargin=!]
\item The weighted walk sum of $G$ satisfies \begin{equation}\label{eq:weightedwalksum}
w(G) \geq \frac{1}{2}\log n - \frac{1}{2}\log(\log_{d - 1}n) - c.
\end{equation}
\item The weighted cycle sum of $G$ satisfies
\begin{equation}\label{eq:weightedcyclesum}
\ell(G) \geq \frac{1}{4}\log n - \frac{1}{2}\log(\log_{d - 1}n) - c.
\end{equation}
\end{enumerate}
\end{theorem}

The coefficient $1/4$ in \eqref{eq:weightedcyclesum} is sharp for each fixed $d$ (Section~\ref{subsec:cycle-sharpness}). Theorem~\ref{thm:walks} may be of independent interest, and immediately implies, via Theorem \ref{thm:main},
that for some constants $a, b > 0$, if $d \geq 4$ is even and $G$ is a simple $n$-vertex $d$-regular graph, then 
\begin{equation}\label{eq:polycorrection}
\varepsilon(G) \geq b\Bigl(\frac{n}{(\log n)^2}\Bigr)^a \cdot S_d(n).
\end{equation}
Thus \eqref{eq:polycorrection} gives a polynomial correction to Schrijver's bound \eqref{eq:Schrijverboundregular}. Furthermore, in terms of the constant $c$ of Theorem~\ref{thm:walks}, one may take $a=\lambda/4=1/4-O(1/d)$ and $b=\min\{1,e^{-c\lambda}(\log(d-1))^{\lambda/2}\}$, which equals $1$ for all sufficiently large $d$. The proof of Theorem~\ref{thm:eo-conditioned-girth} gives the following finite bound.

\begin{corollary}\label{cor:finite-girth}
For every even integer $d\ge4$ and every integer $n\ge(d-1)^{64}$, there exists an $n$-vertex $d$-regular pseudograph $G$ such that
\[
 \varepsilon(G)
 \le \frac{\sqrt{\pi d\,e^{-\gamma}n}}{\sqrt{\log_{d-1}n}-2}\cdot S_d(n)
 \le 2\sqrt{\frac{dn}{\log_{d-1}n}}\cdot S_d(n),
\]
where $\gamma$ is Euler's constant.
\end{corollary}

 It seems plausible that the upper bound in Corollary~\ref{cor:finite-girth} from random $d$-regular graphs conditioned on large girth has the correct order for each fixed even $d\ge4$. We make the following conjecture.
\stepcounter{theorem}
\begin{conjecture}\label{conj:lowerbound}
Fix an even integer $d \geq 4$, and let $G$ be any $n$-vertex $d$-regular pseudograph. Then as $n \rightarrow \infty$, 
\[
\varepsilon(G) \geq (1 - o(1))  \cdot e^{-\gamma/2} \Bigl(\frac{\pi n}{2\log_{d-1} n}\Bigr)^{1/2} \cdot S_d(n).
\]
\end{conjecture}
The case $d=n-1$ lies outside the range of Corollary~\ref{cor:finite-girth}. In this case, McKay~\cite{McK90} proved that for each fixed $\eta>0$, as $n$ tends to infinity through odd integers,
\[
  \varepsilon(K_n)=\bigl(1+O(n^{-1/2+\eta})\bigr)\sqrt{n/e}
   \left(\frac{2^{n+1}}{\pi n}\right)^{(n-1)/2}.
\]
The calculation in Section~\ref{subsec:complete-graph} shows that the lower bound in \eqref{eq:cyclebound} is  
$$(1 + o(1))\sqrt{\frac{2}{\pi}}\,(2e^\gamma)^{1/4}\,n^{-3/4} \cdot 
\varepsilon(K_n).$$

\subsection{Upper bounds.} 
Theorem \ref{thm:main} suggests that 
the number of Eulerian orientations in a $d$-regular graph is higher when there are many short cycles. For even $d$, it is not hard to see that the maximum value of $\varepsilon(G)$ when $G$ is an $n$-vertex $d$-regular pseudograph 
is $2^{dn/2}$, with equality for the graph obtained from $n$ isolated vertices by adding $d/2$ loops at each vertex. If we restrict to $d$-regular multigraphs and $n$ is even, it is straightforward to show that $\varepsilon(G)$ is maximized when every component consists of $d$ repeated edges between two vertices. The problem becomes more challenging when restricting to graphs (no loops or multiple edges). For upper bounds on regular graphs see also Las Vergnas~\cite{LV90}. Schrijver~\cite[Eq.~(22)]{Sch83} proposed the following upper bounds for simple graphs.

\needspace{5\baselineskip}

\begin{conjecture} {\bf (Schrijver~\cite{Sch83})} \label{conj:schrijver-degrees}
Let $G$ be an $n$-vertex graph with even degrees $d_1,d_2,\dots,d_n$, with no loops or multiple edges. Then
\[
\varepsilon(G)\le \prod_{i=1}^n \varepsilon(K_{d_i+1})^{1/(d_i+1)}.
\]
\end{conjecture}
For $d$-regular graphs, Conjecture~\ref{conj:schrijver-degrees} is easily proved for $d = 2$, where a disjoint union of triangles provides a tight example. For $d=4$, the conjecture and its equality case follow from the bounds and enumeration of Bartzos and Samaris~\cite[Theorems~2--3 and Table~1]{BS26}. Their bounds are below $24^{n/5}$ for connected graphs on $n\ge16$ vertices, and their table handles $5\le n\le15$. Isaev, McKay and Zhang~\cite{IMZ} had verified the conjecture through 19 vertices. We prove the following more general theorem. In this theorem, we write $K_r^{(2)}$ for the graph obtained by doubling every edge of $K_r$.  

\begin{theorem}\label{thm:main2}
  Let $G$ be an $n$-vertex multigraph with exactly $k$ vertices of degree $2$, exactly $n-k$ vertices of degree $4$, and exactly $m_i$ unordered pairs of distinct vertices joined by precisely $i$ edges, for $i=2,3,4$. Then
$$\varepsilon(G) \leq \varepsilon(K_5)^{\frac{n}{5} - \frac{m_2}{5} - \frac{2m_3}{5} - \frac{2m_4}{5} - \frac{k}{10}}\cdot \varepsilon(K_3^{(2)})^{\frac{m_2}{3} - \frac{k}{6}} \cdot \varepsilon(K_2^{(2)})^{\frac{k}{2}} \cdot 6^{m_3 + m_4}.$$
\end{theorem}

Since $\varepsilon(K_3^{(2)}) = 10$ and $\varepsilon(K_2^{(2)}) = 2$ and 
$\varepsilon(K_5)^{n/5} = 24^{n/5}$, the upper bound in Theorem \ref{thm:main2} for multigraphs with $m=m_2$ double edges, $m_3=m_4=0$, and $k$ vertices of degree two is exactly 
\begin{equation}\label{eq:ebound}\varepsilon(G) \leq  24^{\frac{n}{5} - \frac{m}{5} - \frac{k}{10}} \cdot 10^{\frac{m}{3} - \frac{k}{6}} \cdot 2^{\frac{k}{2}}.
\end{equation}
Taking $k = m = 0$ in \eqref{eq:ebound}, we obtain Schrijver's Conjecture for $d = 4$. The proof of Theorem \ref{thm:main2} also shows that the only examples which achieve equality in Schrijver's Conjecture for $d = 4$ are graphs whose components are all isomorphic to $K_5$. More generally, the graph $G$ consisting of $m/3 - k/6$ copies of $K_3^{(2)}$, $k/2$ copies of $K_2^{(2)}$, and $n/5 - m/5 - k/10$ copies of $K_5$ has exactly as many Eulerian orientations as the upper bound in Theorem \ref{thm:main2}, whenever the three component counts are nonnegative integers.

It would be very interesting to determine whether there exists an upper bound on $\varepsilon(G)$ which is of a similar form to the lower bound given by Theorem \ref{thm:main}. 
\begin{conjecture}\label{conj:cycle-upper}
For each fixed even $d\ge4$ there exists a constant $c_d$ such that for any $n$-vertex $d$-regular pseudograph $G$,
$$\varepsilon(G) \leq \exp(c_d \cdot \ell(G))  \cdot S_d(n).$$
\end{conjecture}
The constant must depend on $d$ in this class. For two vertices joined by $d$ parallel edges, $\ell(G)=d/(2(d-1))$, whereas $\varepsilon(G)=\bigl(2^d/\binom d{d/2}\bigr)\cdot S_d(2)\sim\sqrt{\pi d/2}\cdot S_d(2)$.

\subsection[Explicit examples]{Explicit examples.}
Let $G=C_L\mathbin{\square}\cdots\mathbin{\square} C_L$ be the Cartesian
product of $k$ cycles of length $L$, so $n=L^k$ and $d=2k$.
For $k=2$, Lieb~\cite{Lie67} proved that
$\lim_{L\to\infty}\varepsilon(G)^{1/n}=(4/3)^{3/2}=1.53960\ldots$.
The lower bound for this limit from Theorem~\ref{thm:main} is approximately
$1.52238$. For $k=3$, Isaev, McKay and Zhang~\cite[Section~6.2]{IMZ}
estimate $\lim_{L\to\infty}n^{-1}\log\varepsilon(G)$ as
$0.9252\pm0.0002$, giving approximately $2.5224$ for
$\lim_{L\to\infty}\varepsilon(G)^{1/n}$. We estimate the lower bound from
Theorem~\ref{thm:main} at $2.51675$. For $k=4$, the value of
$\lim_{L\to\infty}\varepsilon(G)^{1/n}$ is unknown, while the lower bound
from Theorem~\ref{thm:main} is approximately $4.39007$.

To evaluate the weighted cycle sum per vertex $\ell(G)/n$ from
\eqref{eq:cyclesum}, we use the polygon counts tabulated by Clisby, Liang
and Slade~\cite{CLS07}, including Jensen's square-lattice counts. Their
tables give polygon counts through lengths $110$, $32$ and $26$ in
dimensions $2$, $3$ and $4$, respectively.

For a triangular lattice $T$ on the torus, let $n=|V(T)|$ and let $m$
be the length of its shortest non-contractible cycle. Baxter's
solution~\cite{Baxter,LW72} gives, as $m\to\infty$,
\[
 \varepsilon(T)=\left(\frac{3\sqrt{3}}{2}\right)^{(1+o_m(1))n},
 \qquad \frac{3\sqrt{3}}{2}=2.598076\ldots .
\]
For all sufficiently large $m$, every ball of a fixed radius agrees with the
corresponding triangular-lattice ball. Hence $\varepsilon(T)^{1/n}$ has
the same limit along every sequence of such tori with $m\to\infty$ by Theorem~1.3 of Bencs, Borb\'enyi and Csikv\'ari~\cite{BBC24}.

\begin{center}
\begin{tikzpicture}[
  x={(.58cm,0cm)}, y={(.29cm,.502295cm)},
  lattice edge/.style={draw=black!64,line width=.38pt},
  outer edge/.style={draw=black!38,line width=.32pt},
  ball edge/.style={draw=teal!58!black,line width=.85pt},
  seam/.style={draw=black!84,line width=.72pt},
  gluing arrow/.style={->,>=stealth,draw=black,line width=.95pt}]
  \def\ballboundary{(4,0)--(0,4)--(-4,4)--(-4,0)--(0,-4)--(4,-4)--cycle}
  \def\drawlatticeball{%
    \fill[teal!10] \ballboundary;
    \foreach \i in {-4,...,4}{
      \pgfmathtruncatemacro{\lowerend}{max(-4,-4-\i)}
      \pgfmathtruncatemacro{\upperend}{min(4,4-\i)}
      \draw[lattice edge] (\i,\lowerend)--(\i,\upperend);
      \draw[lattice edge] (\lowerend,\i)--(\upperend,\i);
      \draw[lattice edge] ({\i+\lowerend},{-\lowerend})
                         --({\i+\upperend},{-\upperend});
    }
    \draw[ball edge] \ballboundary;
    \foreach \i in {-4,...,4}{
      \pgfmathtruncatemacro{\lowerend}{max(-4,-4-\i)}
      \pgfmathtruncatemacro{\upperend}{min(4,4-\i)}
      \foreach \j in {\lowerend,...,\upperend}
        \fill[black!85] (\i,\j) circle[radius=.82pt];
    }
    \fill[teal!55!black] (0,0) circle[radius=2.1pt];
  }
  \drawlatticeball
  \node[anchor=north,font=\small] at (2.7,-5.4) {Infinite triangular lattice};
  \begin{scope}[shift={(13.3,0)}]
    \begin{scope}[x=1.13cm,y=1.13cm,line cap=round,line join=round]
      \def\paintsurface{%
        \foreach \angle in {-90,-88,...,90}{
          \pgfmathsetmacro{\height}{.77*sin(\angle)*cos(37)}
          \pgfmathsetmacro{\outerradius}{2.5+.77*cos(\angle)}
          \pgfmathsetmacro{\innerradius}{2.5-.77*cos(\angle)}
          \pgfmathsetmacro{\lightness}{50+50*sin(\angle)}
          \path[fill=top tone!\lightness!bottom tone,even odd rule]
            (0,\height) ellipse[x radius=\outerradius,y radius={\outerradius*sin(37)}]
            (0,\height) ellipse[x radius=\innerradius,y radius={\innerradius*sin(37)}];
        }}
      \colorlet{bottom tone}{black!32}
      \colorlet{top tone}{black!1}
      \paintsurface
      \begin{scope}[even odd rule]
        \clip
      (-0.246,-2.268)--(0.239,-2.272)--(0.246,-2.265)--(0.357,-2.265)--(0.365,-2.257)--
      (0.446,-2.257)--(0.453,-2.25)--(0.646,-2.235)--(1.119,-2.146)--(1.519,-2.028)--
      (1.815,-1.909)--(2.067,-1.784)--(2.252,-1.673)--(2.54,-1.458)--(2.766,-1.24)--
      (2.944,-1.018)--(3.099,-0.751)--(3.151,-0.633)--(3.21,-0.455)--(3.247,-0.293)--
      (3.262,-0.196)--(3.262,-0.122)--(3.269,-0.115)--(3.262,0.196)--(3.232,0.366)--
      (3.173,0.573)--(3.099,0.751)--(2.988,0.951)--(2.884,1.099)--(2.792,1.191)--
      (2.781,1.187)--(2.825,1.106)--(2.855,1.017)--(2.87,0.906)--(2.855,0.751)--
      (2.81,0.61)--(2.714,0.44)--(2.581,0.277)--(2.407,0.118)--(2.207,-0.03)--
      (1.933,-0.193)--(1.63,-0.341)--(1.282,-0.481)--(0.912,-0.607)--(0.683,-0.666)--
      (0.365,-0.718)--(0.194,-0.725)--(0.187,-0.733)--(-0.19,-0.733)--(-0.198,-0.725)--
      (-0.427,-0.711)--(-0.686,-0.666)--(-0.886,-0.614)--(-1.064,-0.555)--(-1.204,-0.496)--
      (-1.36,-0.415)--(-1.53,-0.296)--(-1.659,-0.167)--(-1.704,-0.1)--(-1.733,-0.026)--
      (-1.733,0.033)--(-1.689,0.174)--(-1.593,0.344)--(-1.441,0.51)--(-1.241,0.658)--
      (-1.012,0.777)--(-0.768,0.866)--(-0.523,0.925)--(-0.209,0.965)--(-0.216,1.343)--
      (-0.224,1.35)--(-0.224,1.461)--(-0.231,1.469)--(-0.253,1.787)--(-0.305,2.127)--
      (-0.335,2.231)--(-0.361,2.264)--(-0.753,2.22)--(-1.152,2.138)--(-1.5,2.035)--
      (-1.767,1.931)--(-2.04,1.798)--(-2.329,1.62)--(-2.558,1.443)--(-2.806,1.195)--
      (-2.917,1.054)--(-3.006,0.921)--(-3.132,0.677)--(-3.206,0.47)--(-3.235,0.351)--
      (-3.258,0.225)--(-3.265,0.055)--(-3.272,0.048)--(-3.258,-0.226)--(-3.198,-0.492)--
      (-3.132,-0.677)--(-3.065,-0.818)--(-2.895,-1.084)--(-2.725,-1.284)--(-2.566,-1.436)--
      (-2.44,-1.539)--(-2.144,-1.739)--(-1.781,-1.924)--(-1.559,-2.013)--(-1.337,-2.087)--
      (-1.056,-2.161)--(-0.642,-2.235)--cycle;
        \colorlet{bottom tone}{teal!45!gray}
        \colorlet{top tone}{teal!8}
        \paintsurface
      \end{scope}
      \def\screenX#1#2{(2.5+.77*cos(#2))*cos(#1)}
      \def\screenY#1#2{(2.5+.77*cos(#2))*sin(#1)*sin(37)+.77*sin(#2)*cos(37)}
      \foreach \offset in {0,180}{
        \def\silhouetteangle{\offset+atan(sin(\t)/tan(37))}
        \draw[draw=black!52,line width=.45pt]
          plot[domain=0:360,samples=101,variable=\t,smooth]
          ({\screenX{\t}{\silhouetteangle}},{\screenY{\t}{\silhouetteangle}});
      }
      \def\drawvisiblearcs#1#2{%
        \foreach \ia/\ja/\di/\dj/\first/\last/\samples in {#2}{
          \def\uangle{-119+36*(\ia+\di*\t)}
          \def\vangle{43+36*(\ja+\dj*\t)}
          \draw[#1] plot[domain=\first:\last,samples=\samples,variable=\t,smooth]
            ({\screenX{\uangle}{\vangle}},{\screenY{\uangle}{\vangle}});
        }}
      \drawvisiblearcs{draw=black!69,line width=.41pt}{
        0/0/0/1/0/2.43743/23,0/0/0/1/7.43743/10/25,
        1/0/0/1/0/2.33906/23,1/0/0/1/7.33906/10/25,
        2/0/0/1/0/2.57935/25,2/0/0/1/7.57935/10/23,
        3/0/0/1/0/3.41085/32,3/0/0/1/8.41085/10/16,
        4/0/0/1/0/4.61903/43,4/0/0/1/9.61903/10/5,
        5/0/0/1/0.17368/5.17368/47,6/0/0/1/0.27205/5.27205/46,
        7/0/0/1/0.03176/5.03176/46,8/0/0/1/0/4.20026/39,
        8/0/0/1/9.20026/10/9,9/0/0/1/0/2.99208/28,
        9/0/0/1/7.99208/10/20,0/0/1/0/0/4.54567/42,
        0/0/1/0/7.06544/10/28,0/1/1/0/0/10/91,
        0/2/1/0/0/10/91,0/3/1/0/2.61973/8.99138/59,
        0/4/1/0/3.45302/8.15809/44,0/5/1/0/4.54567/7.06544/24,
        0/8/1/0/0/2.61973/25,0/8/1/0/8.99138/10/11,
        0/9/1/0/0/3.45302/33,0/9/1/0/8.15809/10/18,
        0/0/1/-1/0/2.27123/22,0/0/1/-1/4.86672/10/48,
        0/1/1/-1/0/2.80953/27,0/1/1/-1/5.72325/10/40,
        0/2/1/-1/0/3.25779/31,0/2/1/-1/6.94229/8.73079/18,
        0/2/1/-1/9.06913/10/10,0/3/1/-1/0.66362/3.69702/29,
        0/4/1/-1/1.57277/4.2085/25,0/5/1/-1/2.27123/4.86672/25,
        0/6/1/-1/2.80953/5.72325/28,0/7/1/-1/3.25779/6.94229/35,
        0/7/1/-1/8.73079/9.06913/5,0/8/1/-1/0/0.66362/7,
        0/8/1/-1/3.69702/10/58,0/9/1/-1/0/1.57277/16,
        0/9/1/-1/4.2085/10/54}
      \drawvisiblearcs{ball edge}{
        4/0/-1/1/0/2.42723/23,0/4/-1/0/1.84191/4/21,
        -4/4/0/-1/0/3.72795/35,4/-4/0/1/3.61903/4/5}
      \foreach \i/\j in {
        0/0,0/1,0/2,0/8,0/9,1/0,1/1,1/2,1/8,1/9,2/0,2/1,
        2/2,2/8,2/9,3/0,3/1,3/2,3/3,3/9,4/0,4/1,4/2,4/3,
        4/4,5/1,5/2,5/3,5/4,5/5,6/1,6/2,6/3,6/4,6/5,7/1,
        7/2,7/3,7/4,7/5,8/0,8/1,8/2,8/3,8/4,9/0,9/1,9/2,
        9/8,9/9}
        \fill[black!83] ({\screenX{-119+36*\i}{43+36*\j}},
                          {\screenY{-119+36*\i}{43+36*\j}}) circle[radius=.76pt];
      \fill[teal!55!black] ({\screenX{-119}{43}},{\screenY{-119}{43}})
        circle[radius=2.1pt];
    \end{scope}
    \node[anchor=north,font=\small] at (2.7,-5.4) {$10\times10$ triangular torus};
  \end{scope}
\end{tikzpicture}

{\small The shaded regions are the induced balls of radius $4$ about the
marked vertices.\par}
\end{center}

We now compare this with the lower bound given by Theorem~\ref{thm:main}.
Using Jensen's enumeration~\cite[Table~1]{Jensen04} gives, for any vertex
$x\in V(T)$ and $m>60$,
\[
 \frac{1}{n}\log\varepsilon(T)
 \ge \log(5/2)+5\log(6/5)\sum_{i=3}^{60}\frac{\ell_i(T,x)}{i\,5^i}
 =0.93918455\ldots .
\]
Here $\ell_i(T,x)$ counts the $i$-cycles containing $x$. For $i<m$, these
cycles lift bijectively to the $i$-cycles through a fixed lift of $x$ in the
triangular lattice. By vertex transitivity, $n\ell_i(T,x)=i\ell_i(T)$, since each $i$-cycle
contains $i$ vertices.

For the regular hyperbolic tiling $\{p,d\}$, each face has $p$ sides and $d$
faces meet at each vertex, where $1/p+1/d<1/2$ and $d\ge4$ is even.
Let $G$ be an $n$-vertex quotient on a closed orientable surface, and let $m$
be the length of its shortest non-contractible cycle.
\v{S}ir\'a\v{n}~\cite{Siran01} constructed
such quotients with arbitrarily large $m$. For fixed $p$ and $d$, these
quotients converge locally to the tiling as $m\to\infty$, so
$\varepsilon(G)^{1/n}$ has the same limit along every such sequence
by Bencs, Borb\'enyi and Csikv\'ari~\cite[Theorem~1.3]{BBC24}.

\begin{center}
\begin{tikzpicture}[scale=2.44,line cap=round,line join=round]
  \def\twofaces{(.321969,-.233924)
    arc[start angle=99,delta angle=18,radius=1.495349]
    arc[start angle=27,delta angle=18,radius=1.495349]
    arc[start angle=315,delta angle=18,radius=1.495349]
    arc[start angle=243,delta angle=18,radius=1.495349]
    arc[start angle=261,delta angle=13.717474,radius=1.495349]
    arc[start angle=184.717474,delta angle=40.282526,radius=.353004]
    arc[start angle=135,delta angle=40.282526,radius=.353004]
    arc[start angle=85.282526,delta angle=13.717474,radius=1.495349]
    --cycle}
  \begin{scope}
    \clip (0,0) circle[radius=.95];
    \fill[teal!10] \twofaces;
    \foreach \turn in {0,72,144,216,288}{
      \foreach \angle/\dist/\rad in {
        0/1.798908/1.495349, 13.61382/1.060477/.353004, -13.61382/1.060477/.353004,
        36/1.019828/.200123, 15.51011/1.006899/.117668, -15.51011/1.006899/.117668,
        6.64537/1.005155/.101668, -6.64537/1.005155/.101668,
        23.73361/1.003466/.083333, -23.73361/1.003466/.083333,
        33.39097/1.002079/.064509, -33.39097/1.002079/.064509, 0/1.001601/.056607}{
        \draw[black!64,line width=.38pt]
          ({\turn+\angle}:\dist) circle[radius=\rad];
      }
    }
    \draw[teal!58!black,line width=.85pt] \twofaces;
    \foreach \vx/\vy in {0.321969/-0.233924,-0.122981/-0.378497,-0.397975/0.000000,-0.122981/0.378497,
                       0.321969/0.233924,0.678874/0.220580,0.781071/0.000000,0.678874/-0.220579}
      \fill[teal!58!black] (\vx,\vy) circle[radius=.82pt];
  \end{scope}
  \draw[black!38,line width=.32pt] (0,0) circle[radius=.95];
\end{tikzpicture}

{\small The $\{5,4\}$ tiling in the Poincar\'e disk.
The boundary of the two shaded pentagons is an $8$-cycle.\par}
\end{center}

Suppose $m>2p-2$. The $dn/p$ faces give distinct $p$-cycles. Each of the
$dn/2$ edges also gives a $(2p-2)$-cycle bounding the two faces that share
it. These boundaries are simple in the quotient since their lengths are
less than $m$, and each determines its pair of faces upon lifting to the
tiling. Keeping these two terms in Theorem~\ref{thm:main} gives
\[
 \varepsilon(G)^{1/n}\ge
 \exp\!\left((d-1)\log\!\frac{d}{d-1}
 \left(\frac{d}{p(d-1)^p}+\frac{d}{2(d-1)^{2p-2}}\right)\right)
 \cdot S_d(n)^{1/n}.
\]
Thus $\varepsilon(G)^{1/n}$ exceeds $1.50466$ for $\{5,4\}$, $2.50591$
for $\{4,6\}$, and $4.41377$ for $\{3,8\}$. In $\{3,8\}$, three consecutive
triangles around a vertex bound a pentagon. The pentagon determines their
common vertex, and there are eight such triples at each vertex. When $m>5$,
these $8n$ pentagons contribute a further $8/7^5$ to $\ell(G)/n$, giving
$\varepsilon(G)^{1/n}>4.41574$. For a vertex $x$ of the infinite tiling
$\mathcal T=\{p,d\}$, Panagiotis~\cite[Theorem~1]{Panagiotis24}
proved that $\limsup_{i\to\infty}\ell_i(\mathcal T,x)^{1/i}<d-1$.
Consequently the tail of the weighted cycle sum
$\sum_{i\ge p}\ell_i(\mathcal T,x)/(i(d-1)^i)$ decays exponentially.

For fixed $d$, one can also compare the bound with the actual count as
$p\to\infty$. The tiling has girth $p$, so its quotient has girth $p$
when $m>2p-2$. Each circuit in an edge partition then uses at least $p$
of the $dn/2$ edges. By the circuit-partition identity~\eqref{eq:circuit-partition-mean},
\[
 \varepsilon(G)^{1/n}\le 2^{d/(2p)}\cdot S_d(n)^{1/n}.
\]
The per-vertex lower bound in Theorem~\ref{thm:main} is at least
$S_d(n)^{1/n}$, so the ratio of the actual per-vertex count to that bound
lies between $1$ and $2^{d/(2p)}$ and tends to $1$.

The diamond lattice is the $4$-regular oxygen network of cubic ice.
Let $G$ be an $n$-vertex periodic quotient whose shortest non-contractible
cycle has length greater than $34$. For any vertex $x$ and $i\le34$,
$\ell_i(G,x)/i$ is the number of $i$-cycles per vertex in the infinite
diamond lattice. Clisby's counts in
Table~16.12 of Guttmann and Jensen~\cite{GJ09}
give $\sum_{i=6}^{34}\ell_i(G,x)/(i\,3^i)=0.00414409148984\ldots$.
Hence Theorem~\ref{thm:main} gives
\[
 \varepsilon(G)^{1/n}\ge
 \exp\!\left(3\log(4/3)\sum_{i=6}^{34}\frac{\ell_i(G,x)}{i\,3^i}\right)
 \cdot\frac32
 =1.5053744189\ldots .
\]
Kolafa's Monte Carlo estimate~\cite{Kolafa14}
for $\lim_{m\to\infty}\varepsilon(G)^{1/n}$, with $m$ the length of the
shortest non-contractible cycle, is $1.5074660(36)$. The bound above is
within $0.14\%$ of this estimate and accounts for about $72\%$ of its excess
over $3/2$.

\section{Proof of Theorem~\ref{thm:main}}
\label{sec:cycle-induction}

\subsection{A lower bound from closed trails.}
At a vertex of even degree, a pairing is a partition of its incidences into
unordered pairs. A loop supplies two distinct incidences. Write $\mathcal P(G)$
for the set of choices of one pairing at every vertex. Given $P\in\mathcal P(G)$,
follow each edge on arrival by the edge whose incidence is paired with the arrival
incidence. This partitions the edges into closed trails. Each trail retains its
cyclic order, and we identify its two directions of traversal. Let $X(P)$ be the
number of trails in this partition.

The circuit-partition identity of Lieb and Wu~\cite[Eq.~(11)]{LW72} and
Las Vergnas~\cite{LV83}, also given by Isaev, McKay and
Zhang~\cite[Theorem~5.1]{IMZ}, states that for a pseudograph with even degrees,
\begin{equation}\label{eq:circuit-partition}
  \sum_{P\in\mathcal P(G)}2^{X(P)}
   =\varepsilon(G)\prod_{v\in V(G)}\bigl(d(v)/2\bigr)! .
\end{equation}
Indeed, both sides count pairs $(P,O)$ in which $O$ is an Eulerian orientation
that directs each trail of $P$ consistently. For a fixed $P$, each of its $X(P)$
trails can be directed independently in either of its two directions. For a
fixed $O$, a compatible pairing at $v$ is a bijection from the $d(v)/2$ incoming
incidences to the $d(v)/2$ outgoing incidences. There are $(d(v)/2)!$ such
bijections at $v$, and the choices at different vertices are independent.

Now suppose $G$ is $d$-regular on $n$ vertices, with $d\ge2$ even. Choose
$P$ uniformly from $\mathcal P(G)$. There are $(d-1)!!$ pairings at each vertex,
so $|\mathcal P(G)|=((d-1)!!)^n$. Dividing \eqref{eq:circuit-partition}
by this number and using
\[
 \frac{(d-1)!!}{(d/2)!}
 =\frac{d!}{2^{d/2}((d/2)!)^2}
 =2^{-d/2}\binom{d}{d/2}
\]
gives the form of the identity used below,
\begin{equation}\label{eq:circuit-partition-mean}
  \varepsilon(G)=\EE\!\left(2^{X(P)}\right)\cdot S_d(n).
\end{equation}

Let $w_i^\circ(G)$ count the rooted, directed closed trails of length $i$ in $G$.
Such a trail starts with a specified directed edge, returns to its starting vertex
and uses each underlying edge copy at most once. Thus its length is at most
$|E(G)|=dn/2$. Define the {\em weighted trail sum}
\[
 w^\circ(G) :=\sum_{i=1}^{dn/2}\frac{w_i^\circ(G)}{2i(d-1)^i}.
\]

\begin{theorem}\label{thm:trail-bound}
Let $d\ge2$ be even and let $G$ be a $d$-regular pseudograph on $n\ge1$
vertices. Then $\varepsilon(G)\ge 2^{w^\circ(G)}\cdot S_d(n)$, with equality
for $d=2$ and strict inequality for $d\ge4$.
\end{theorem}

\begin{proof}
Let $\mathcal P(G)$ be the set of all tuples $(P_v)_{v\in V(G)}$ such that
$P_v$ partitions the $d$ edge incidences at $v$ into unordered pairs.
A loop at $v$ contributes two distinct incidences. Choose
$P=(P_v)_{v\in V(G)}$ uniformly from $\mathcal P(G)$. Thus the $P_v$
are independent, and each is uniform over the pairings at $v$.
On arrival at $v$, follow the edge incidence paired with the arrival
incidence by $P_v$. Write $X(P)$ for the number of closed trails in the
resulting partition of $E(G)$.

An unrooted, unoriented closed trail $T$ of length $i$ prescribes $r_v$
disjoint pairs at vertex $v$, where $r_v$ is its number of visits to $v$
and $\sum_v r_v=i$. To obtain these pairs in a uniform pairing at $v$,
the first specified incidence must be matched to one particular partner
among $d-1$, the next among $d-3$, and so on. Since the vertices are
paired independently,
\[
 \PP(T\text{ is a trail of }P)
 =\prod_{v\in V(G)}\prod_{j=0}^{r_v-1}\frac{1}{d-2j-1}
 \ge (d-1)^{-i}.
\]
Each $T$ has $2i$ rooted, directed representations, one for each choice of
first edge and direction. Summing the indicators over all unrooted, unoriented
closed trails counts each trail in the partition exactly once. Consequently
\[
  \EE(X(P))=\sum_T\PP(T\text{ is a trail of }P)
  \ge\sum_{i=1}^{dn/2}\frac{w_i^\circ(G)}{2i}(d-1)^{-i}
  =w^\circ(G).
\]
We now use the circuit-partition identity \eqref{eq:circuit-partition-mean}
in the first equality below. The first inequality is Jensen's inequality
for $x\mapsto 2^x$. The second follows from the preceding estimate.
\[
 \varepsilon(G)
   =\EE\!\left(2^{X(P)}\right)\cdot S_d(n)
   \ge 2^{\EE(X(P))}\cdot S_d(n)
   \ge 2^{w^\circ(G)}\cdot S_d(n).
\]

For $d=2$ each component is a cycle, so $X(P)=w^\circ(G)$ is the number of
components. For $d\ge4$, decompose the edges into simple cycles and take
the pairings that follow these cycles. At any vertex at least two distinct
cycles meet. Exchanging the partners between one pair from each of these
cycles joins them into a single trail and leaves every other trail unchanged.
These two choices of $P$ have positive probability and give different values
of $X(P)$. Strict convexity of $2^x$ makes Jensen's inequality strict.
\end{proof}

\subsection{Proof of Theorem~\ref{thm:main}.}
Fix an even integer $d\ge2$. By the definition \eqref{eq:cyclesum}, for a
$d$-regular pseudograph $G$ on $n$ vertices we have
\[
  (d-1)\ell(G)=\sum_{i=1}^n\frac{\ell_i(G)}{(d-1)^{i-1}}.
\]
Thus it suffices to prove
\[
  \varepsilon(G)\ge
     \left(1+\frac{1}{d-1}\right)^{(d-1)\ell(G)}\cdot S_d(n).
\]
We induct on $n$. The empty graph has one orientation and no cycles.
Since $S_d(0)=1$, the bound holds with equality when $n=0$. Let $n\ge1$,
and assume the bound holds for every $d$-regular pseudograph on $n-1$ vertices.

Fix a vertex $v$ of $G$, and let $k$ be the number of loops at $v$.
These loops use $2k$ of its $d$ incidences, leaving $d-2k$ incidences
of edges from $v$ to other vertices. Let $\mathcal P_v$ be the set of
pairings of these $d-2k$ incidences. Ordering the incidences and grouping
successive entries into pairs represents each pairing
$2^{d/2-k}(d/2-k)!$ times, so
\[
 |\mathcal P_v|=\frac{(d-2k)!}{2^{d/2-k}(d/2-k)!}.
\]
This formula gives one empty pairing when $k=d/2$.

For $P\in\mathcal P_v$, form $G_P$ by deleting $v$ and its incident edges.
Whenever $P$ pairs the incidences of $xv$ and $vy$, insert an edge between
$x$ and $y$. Label this edge by the pair of original edge copies and retain
their incidences at $x$ and $y$ as its two ends. If $x=y$, the new edge is a
loop and the two ends remain distinct. Every incidence removed at a vertex
other than $v$ is restored in this operation, so $G_P$ is a $d$-regular
pseudograph on $n-1$ vertices. We count each pairing separately even if
several give isomorphic pseudographs.

Call an Eulerian orientation of $G$ compatible with $P$ if each pair
at $v$ consists of one incoming and one outgoing incidence. First fix an
Eulerian orientation of $G$ and count its compatible pairings. Each of the $k$
loops contributes one incoming and one outgoing incidence, so among
the non-loop incidences exactly $d/2-k$ are incoming and $d/2-k$ are
outgoing. Each bijection from the former to the latter gives one
compatible pairing. Hence each orientation is compatible with
$(d/2-k)!$ pairings.

Conversely, fix $P$ and an Eulerian orientation of $G_P$. If a new edge
is directed from its end at $x$ to its end at $y$, direct its original
edges from $x$ to $v$ and from $v$ to $y$. This rule also applies when
$x=y$, since the direction of the new loop specifies which of its two
ends is the tail. Every pair now contributes one edge in and one edge
out at $v$, and all other vertices keep their incoming and outgoing
incidences. Each of the $k$ deleted loops at $v$ can be directed in two
ways. These are exactly the $2^k$ compatible orientations of $G$
inducing the chosen orientation of $G_P$. Counting compatible pairs
of an orientation and a pairing in these two ways gives
\[
 (d/2-k)!\,\varepsilon(G)=2^k\sum_{P\in\mathcal P_v}\varepsilon(G_P).
\]
From now on $\EE_P$ denotes expectation over a uniform $P\in\mathcal P_v$.
Substituting the value of $|\mathcal P_v|$ into the last equality gives
\begin{equation}\label{eq:cycle-smoothing}
 \varepsilon(G)=\frac{2^k|\mathcal P_v|}{(d/2-k)!}\,
                  \EE_P(\varepsilon(G_P))
       =2^k\,2^{-(d/2-k)}\binom{d-2k}{d/2-k}\,
                  \EE_P(\varepsilon(G_P)) .
\end{equation}

We next compare $\ell(G_P)$ with $\ell(G)$. Every cycle of $G$ avoiding
$v$ remains a cycle of the same length in $G_P$. Consider a cycle $C$
of length $i\ge2$ through $v$. It uses two distinct non-loop incidences
at $v$. If $P$ matches these two incidences, the new edge replaces the
two edges of $C$ at $v$, leaving a cycle of length $i-1$. In particular,
a cycle of length two becomes a loop. For a fixed $P$, each shortened
cycle contains exactly one new edge. Its label recovers the two original
edges and hence the original cycle. The other surviving cycles contain
no new edge. Thus all these contributions to $\ell(G_P)$ are distinct.
A cycle avoiding $v$ has length at most $n-1$, and a shortened cycle has
length $i-1\le n-1$, so all are within the cutoff defining $\ell(G_P)$.

Suppose first that $k<d/2$. In a uniform pairing, either of the two
specified incidences of $C$ has $d-2k-1$ equally likely partners.
Therefore $C$ shortens with probability $1/(d-2k-1)$. Write $|C|$
for the length of a cycle. By linearity of expectation,
\begin{align*}
 \EE_P(\ell(G_P))
 &\ge \sum_{C:\,v\notin V(C)}\frac{1}{(d-1)^{|C|}}
   +\sum_{\substack{C:\,v\in V(C)\\ |C|\ge2}}
         \frac{1}{d-2k-1}\frac{1}{(d-1)^{|C|-1}}\\
 &\ge \sum_{C:\,v\notin V(C)}\frac{1}{(d-1)^{|C|}}
   +\sum_{\substack{C:\,v\in V(C)\\ |C|\ge2}}
         \frac{1}{(d-1)^{|C|}}
  =\ell(G)-\frac{k}{d-1}.
\end{align*}
Here the sums run over cycles of $G$. The second inequality uses
$d-2k-1\le d-1$, and the omitted cycles are the $k$ loops at $v$, each
of weight $1/(d-1)$. If $k=d/2$, the vertex $v$ has only loops. The
unique pairing then gives $G_P=G-v$, and precisely these $k$ loops are
lost. We have proved in both cases that
\begin{equation}\label{eq:cycle-weight-loss}
  (d-1)\EE_P(\ell(G_P))\ge (d-1)\ell(G)-k .
\end{equation}

The induction hypothesis applies to each $G_P$ individually and gives
\[
 \varepsilon(G_P)\ge
    \left(1+\frac{1}{d-1}\right)^{(d-1)\ell(G_P)}\cdot S_d(n-1).
\]
Averaging these inequalities over the finite set $\mathcal P_v$ yields
\[
 \EE_P(\varepsilon(G_P))\ge
   \left(\frac1{|\mathcal P_v|}\sum_{P\in\mathcal P_v}
    \left(1+\frac{1}{d-1}\right)^{(d-1)\ell(G_P)}\right)\cdot S_d(n-1).
\]
The function $x\mapsto (1+1/(d-1))^x$ is convex and increasing, since
its base exceeds one. Applying Jensen's inequality to the displayed
average of powers, and then using \eqref{eq:cycle-weight-loss}, gives
\begin{align*}
 \frac1{|\mathcal P_v|}\sum_{P\in\mathcal P_v}
   \left(1+\frac{1}{d-1}\right)^{(d-1)\ell(G_P)}
 &\ge \left(1+\frac{1}{d-1}\right)^{
          \frac1{|\mathcal P_v|}\sum_{P\in\mathcal P_v}(d-1)\ell(G_P)}\\
 &=\left(1+\frac{1}{d-1}\right)^{(d-1)\EE_P(\ell(G_P))}\\
 &\ge\left(1+\frac{1}{d-1}\right)^{(d-1)\ell(G)-k}.
\end{align*}

The factor taking $S_d(n-1)$ to $S_d(n)$ is $2^{-d/2}\binom{d}{d/2}$.
Dividing the coefficient in \eqref{eq:cycle-smoothing} by this factor gives
\begin{align*}
 \frac{2^k\,2^{-(d/2-k)}\binom{d-2k}{d/2-k}}
      {2^{-d/2}\binom{d}{d/2}}
 &=4^k\frac{\binom{d-2k}{d/2-k}}{\binom{d}{d/2}}\\
 &=\prod_{j=0}^{k-1}
   \frac{4(d/2-j)^2}{(d-2j)(d-2j-1)}
  =\prod_{j=0}^{k-1}\left(1+\frac{1}{d-2j-1}\right)\\
 &\ge \left(1+\frac{1}{d-1}\right)^k .
\end{align*}
For $k=0$ the products are empty and the ratio is one. For $k\ge1$,
each denominator $d-2j-1$ lies between $1$ and $d-1$, which proves the
last inequality. Substituting the lower bounds on the coefficient and
on $\EE_P(\varepsilon(G_P))$ into \eqref{eq:cycle-smoothing} yields
\begin{align*}
 \varepsilon(G)
 &\ge 2^{-d/2}\binom{d}{d/2}
     \left(1+\frac{1}{d-1}\right)^k
     \left(1+\frac{1}{d-1}\right)^{(d-1)\ell(G)-k}\cdot S_d(n-1)\\
 &= \left(1+\frac{1}{d-1}\right)^{(d-1)\ell(G)}\cdot S_d(n)\\
 &= \prod_{i=1}^n
      \left(1+\frac{1}{d-1}\right)^{\ell_i(G)/(d-1)^{i-1}}\cdot S_d(n),
\end{align*}
which is \eqref{eq:cyclebound}. This completes the induction.\qed

\section{Linear algebra}
\label{sec:spectral}

\subsection{Non-backtracking matrices.}
Let $H$ be a pseudograph on $N$ vertices with $m$ edges. The adjacency
matrix $A$ counts edge multiplicities off the diagonal and twice the number
of loops on the diagonal. Let $D$ be the diagonal matrix of vertex degrees.
Write $\vec E(H)$ for the $2m$ directed edge copies, obtained by choosing an
edge copy and designating one of its two incidences as its tail. Write
$e^{-1}$ for the other direction of the same edge copy, including when $e$
is a loop. The {\em non-backtracking matrix} (or {\em Hashimoto matrix}~\cite{Has89})
is indexed by $\vec E(H)$ and has entries
\[
 B_{e,f}=\begin{cases}
 1 &\text{if }\head(e)=\tail(f)\text{ and } f\ne e^{-1},\\
 0 &\text{otherwise.}
 \end{cases}
\]
Thus $w_j(H)=\tr(B^j)$, since the trace also checks the transition from
the last directed edge to the first.

We use the Ihara--Bass determinant formula~\cite{Iha66,Bass92}.
For finite graphs with loops and parallel edges, see Kotani and
Sunada~\cite[Theorem~1.1]{KS00}. We also use the Perron--Frobenius
theorem~\cite[Theorem~8.3.1]{HJ13}.

\begin{proposition}\label{prop:nb-spectrum}
Suppose $H$ has maximum degree at most $d$, where $d\ge3$.
Then, as an identity of rational functions,
\[
 \det(tI-B)=(t^2-1)^{m-N}\det(t^2I-tA+D-I).
\]
Every eigenvalue of $B$ has modulus at most $d-1$.
At most $2N$ eigenvalues, counted with algebraic multiplicity, are
nonreal, and each of these has modulus at most $\sqrt{d-1}$.
If $m\ge1$, the spectral radius $\rho(B)$ is a real eigenvalue.
If $H$ has average degree $\bar d>2$, then $\rho(B)\ge\bar d-1$.
\end{proposition}

\begin{proof}
Put $Q(t)=t^2I-tA+D-I$. The cited determinant formula applies componentwise
when the minimum degree is at least two. Adding an isolated vertex leaves
both sides unchanged, since it multiplies $\det Q(t)$ by $t^2-1$.
Adding a leaf multiplies both sides by $t^2$. On the left, expand first
along the row of the directed edge entering the leaf, then along the
column of its reverse. On the right, eliminating the leaf's diagonal
entry $t^2$ cancels the extra one in the degree of its neighbor.
Deleting leaves and isolated vertices therefore proves the identity in
general.

Row $e$ of $B$ sums to $\deg(\head(e))-1\le d-1$, giving the
eigenvalue bound. When $m\ge1$, Perron--Frobenius gives the real eigenvalue
$\rho(B)$. Away from $\pm1$ the eigenvalues of $B$ are the roots of
$\det Q(t)$ with the same multiplicities. This polynomial is monic of
degree $2N$. For a nonreal root $\mu$, choose a unit vector $x$ with
$Q(\mu)x=0$. The scalar quadratic
\[
 t^2-(x^*Ax)t+x^*(D-I)x
\]
has real coefficients and roots $\mu,\bar\mu$. Their product gives
$|\mu|^2=x^*(D-I)x\le d-1$.

Finally, $\mathbf1^{T}Q(t)\mathbf1=N(t-1)(t-\bar d+1)$, so the
smallest eigenvalue of the real symmetric matrix $Q(\bar d-1)$ is
nonpositive. As $Q(t)$ is positive definite for large positive $t$,
continuity gives $\det Q(t)=0$ at some $t\ge\bar d-1$. If $\bar d>2$,
then $t>1$, so the determinant identity gives $\rho(B)\ge t$.
\end{proof}

\subsection{Trace estimates.}
The following lemma pairs consecutive even and odd powers to control
the contribution of negative real eigenvalues.

\begin{lemma}\label{lem:nb-trace}
Let $H$ be a pseudograph on $N$ vertices with at least one edge and
maximum degree at most $d$, where $d\ge3$. Write $q=d-1$, $C=B/q$ and
$c^*=\rho(B)/q$. Let $1\le L_0\le L$, with $L_0$ even and $L$ odd,
and let $a_{L_0}\ge\cdots\ge a_L\ge0$. Then
\[
 \sum_{j=L_0}^L a_j\tr(C^j)\ge
 \sum_{j=L_0}^L a_j(c^*)^j-2N\sum_{j=L_0}^L a_jq^{-j/2}.
\]
\end{lemma}

\begin{proof}
Set $F(x)=\sum_{j=L_0}^L a_jx^j$. For $-1\le x\le1$, each consecutive
pair contributes $x^{2h}(a_{2h}+xa_{2h+1})\ge0$, so $F(x)\ge0$.
The trace of $F(C)$ is the sum of $F(\mu)$ over the eigenvalues $\mu$
of $C$, with algebraic multiplicity. By Proposition~\ref{prop:nb-spectrum}
one term is $F(c^*)$, and all other real eigenvalues lie in $[-1,1]$.
Each nonreal eigenvalue has $|\mu|\le q^{-1/2}$, hence
$|F(\mu)|\le\sum_{j=L_0}^L a_jq^{-j/2}$. There are at most $2N$ such
terms, and their sum is real by conjugacy. This gives the inequality.
\end{proof}

\section{Counting walks and cycles}
\label{sec:walks-cycles}

Let $G$ be a simple $d$-regular graph on $n$ vertices, where $d\ge3$, and
put $L_0=2\lceil\log_{d-1} n\rceil$. Since $d\le n-1$,
\[
  2\le L_0\le 4\log_{d-1} n\le 4\log_2 n.
\]
We prove Theorem~\ref{thm:walks} with $c=10$. For $n<e^{16}$, both
weighted sums are nonnegative and the right sides of \eqref{eq:weightedwalksum}
and \eqref{eq:weightedcyclesum} are at most $\frac12\log n-10<0$.
For the proof of the two bounds, assume $n\ge e^{16}$.

\subsection{A trace estimate.}
Let $H$ be a subgraph of $G$ with at least one edge, and let $B_H$ be
its non-backtracking matrix. Suppose $\rho(B_H)/(d-1)\ge1-\alpha$ for
some $0\le\alpha\le1$. Let $L\ge L_0$ be odd, let
$0\le\beta\le1/(2L)$, and set $b_j=1/(2j)-\beta$. We claim that
\begin{equation}\label{eq:walk-trace-estimate}
  \sum_{j=L_0}^{L} b_j\frac{w_j(H)}{(d-1)^j}
  \ge \frac12\log\frac{L+1}{L_0}-\beta L-\frac{\alpha L}{2}
      -\frac{2}{1-2^{-1/2}}.
\end{equation}
Indeed, $w_j(H)=\tr(B_H^j)$, since the diagonal entries count
non-backtracking walks that return to their initial directed edge.
The coefficients $b_j$ are nonnegative and nonincreasing, so
Lemma~\ref{lem:nb-trace} bounds the left side below by
\[
  \sum_{j=L_0}^{L} b_j(1-\alpha)^j
       -2n\sum_{j\ge L_0}(d-1)^{-j/2}.
\]
The last sum, including its factor $2n$, is at most
$2/(1-2^{-1/2})$ since $(d-1)^{L_0/2}\ge n$.
Now $(1-\alpha)^j\ge1-j\alpha$, $j b_j\le1/2$, and
$\sum_{j=L_0}^L j^{-1}\ge\log((L+1)/L_0)$, proving
\eqref{eq:walk-trace-estimate}.

For $H=G$, all row sums of $B_G$ are $d-1$, so $\rho(B_G)=d-1$.
Taking $\alpha=\beta=0$ and the largest odd $L\le n$ in
\eqref{eq:walk-trace-estimate} and using \eqref{eq:walksum} gives
\[
 w(G)\ge\frac12\log\frac{n}{L_0}-(4+2\sqrt2)
 \ge\frac12\log n-\frac12\log(\log_{d-1}n)-8.
\]
Here $L+1\ge n$, $L_0\le4\log_{d-1}n$ and
$2/(1-2^{-1/2})+\log 2=4+2\sqrt2+\log 2<8$.
This proves \eqref{eq:weightedwalksum}.

\subsection{From walks to cycles.}
For $H\subseteq G$, let $\mathcal C_j(H)$ be the set of simple cycles of length $j$
in $H$, counted up to rotation and reversal, and let $\ell_j(H)$ be their number.
In this subsection a subscript $j$ denotes an exact cycle length, and a subscript
$L$ denotes the cutoff in the weighted sums. For an integer $L\ge1$, put
\[
  \ell_L(H):=\sum_{j=1}^{L}\frac{\ell_j(H)}{(d-1)^j}.
\]
For each $x\in V(H)$, put
\[
  \ell_j(H,x)=\bigl|\{C\in\mathcal C_j(H):x\in V(C)\}\bigr|,
  \qquad
  \ell_L(H,x):=\sum_{j=1}^{L}\frac{\ell_j(H,x)}{(d-1)^j}.
\]
In both sums $d$ is the degree in $G$. If $\ell_L(H,x)\le\delta$ for every
$x\in V(H)$, we claim that
\begin{equation}\label{eq:walks-to-cycles}
  \ell_L(H)\ge\sum_{j=1}^{L}
      \left(\frac{1}{2j}-3\delta\right)\frac{w_j(H)}{(d-1)^j}.
\end{equation}
The $2j\ell_j(H)$ rooted directed traversals of simple $j$-cycles
have total weight $\ell_L(H)$ when summed over $j\le L$.
For any other walk $W=(v_0,\ldots,v_j=v_0)$ counted by $w_j(H)$, read indices modulo $j$ and take
\[
 h=\min\{s\in\{1,\ldots,j-1\}:v_{b+s}=v_b\text{ for some }b\},\qquad v_{a+h}=v_a.
\]
Then $h\ge3$ and $v_a,\ldots,v_{a+h-1}$ are distinct, so $D=(v_a,\ldots,v_{a+h})$
is a simple cycle. After traversing $D$, follow $W$ from $v_{a+h}$ back to
$v_{a+j}=v_a$. This remaining segment has $m=j-h$ directed edges, in order
\[
 e_b=(v_{a+h+b-1},v_{a+h+b})\qquad(1\le b\le m).
\]
Writing $e^{-1}$ for the reverse of $e$, let
\[
 t=\max\{s\in\{0,\ldots,\lfloor m/2\rfloor\}:e_b=e_{m+1-b}^{-1}
                 \text{ for every }1\le b\le s\}.
\]
Set
\[
 P=e_1\cdots e_t \quad \mbox{and} \quad Q=e_{t+1}\cdots e_{m-t}.
\]
Here $P$ is the empty walk if $t=0$. Since $m>0$ and no two consecutive
edges $e_b,e_{b+1}$ are inverses, $m>2t$.
Maximality of $t$ gives $e_{t+1}\ne e_{m-t}^{-1}$. Thus $Q$ is cyclically
non-backtracking and $W$ is a cyclic shift of $DPQP^{-1}$.

To count these walks, first fix the rooted directed walk $Q$, of
length $i$. From its root construct $P$ backwards. There are at most
$p_t=d(d-1)^{t-1}$ choices when $P$ has length $t\ge1$, and $p_0=1$.
At the other end of $P$, the sum of $(d-1)^{-|D|}$ over the choices of
$D$ is at most $\delta$. Each choice has two directions for $D$ and
at most $j$ choices of the original root. These factors cancel the
divisor $2j$. Since $j=i+2t+|D|$, the total weight of walks other
than simple cycles is at most
\[
  \delta\sum_{i=1}^{L}\frac{w_i(H)}{(d-1)^i}
           \sum_{t\ge0}\frac{p_t}{(d-1)^{2t}}
  =\delta\left(1+\frac{d}{(d-1)(d-2)}\right)
           \sum_{i=1}^{L}\frac{w_i(H)}{(d-1)^i}.
\]
The factor in parentheses is at most $3$ for $d\ge3$.
This proves \eqref{eq:walks-to-cycles}.

Let $L$ be the largest odd integer at most $\sqrt n/10$ and
$\sigma=n^{-1/2}$. Since $\sqrt n/\log n$ is increasing for $n\ge e^{16}$,
\[
 L>\frac{\sqrt n}{10}-2\ge\frac{\sqrt n}{20}
   \ge4\log_2 n\ge L_0 .
\] While some vertex $x$ satisfies
$\ell_L(H,x)>\sigma$ in the current graph $H$, delete $x$, starting
with $H=G$. Each deletion removes more than $\sigma$ from $\ell_L(H)$.
If $s$ vertices were deleted and $H$ is the final graph, then
\[
  \ell_L(G)\ge \ell_L(H)+s\sigma,\qquad
  \ell_L(H,x)\le\sigma\quad(x\in V(H)).
\]
We may assume $\ell_L(G)<\log n$, since otherwise
\eqref{eq:weightedcyclesum} follows immediately. Thus
$s/n<(\log n)/\sqrt n\le16e^{-8}<1/4$ and at most $ds$ edges were
removed. The average degree of $H$ is therefore at least
$d-ds/(n-s)>2d/3\ge2$.
By Proposition~\ref{prop:nb-spectrum},
\[
  \frac{\rho(B_H)}{d-1}\ge1-\alpha,\qquad
  \alpha=\frac{ds}{(d-1)(n-s)}.
\]
Here $0\le\alpha\le1$ and $\alpha L/2\le s\sigma$, since
$n-s\ge n/2$, $d/(d-1)\le3/2$, and $L\le\sqrt n/10$.
Apply \eqref{eq:walks-to-cycles} with $\delta=\sigma$, followed by
\eqref{eq:walk-trace-estimate} with $\beta=3\sigma\le1/(2L)$.
The terms with $j<L_0$ have nonnegative coefficients. We obtain
\begin{align*}
  \ell(G)\ge \ell_L(G)
    &\ge s\sigma+\frac12\log\frac{L+1}{L_0}
             -3\sigma L-\frac{\alpha L}{2}-(4+2\sqrt2)\\
    &\ge \frac14\log n-\frac12\log(\log_{d-1} n)-10 .
\end{align*}
Indeed, $s\sigma\ge\alpha L/2$, $3\sigma L\le3/10$, and
$(L+1)/L_0\ge\sqrt n/(80\log_{d-1}n)$. The total additive loss is at most
$4+2\sqrt2+3/10+\frac12\log 80<10$. This proves
\eqref{eq:weightedcyclesum} and completes the proof of
Theorem~\ref{thm:walks}.

\subsection{Sharpness of the leading coefficient.}
\label{subsec:cycle-sharpness}
Fix $d\ge3$ and let $n$ tend to infinity with $dn$ even. Pair the $dn$
labelled incidences of $n$ vertices uniformly, and let $G^*$ be the
resulting pseudograph. Choosing $i$ vertices in cyclic order, two
ordered incidences at each vertex, and requiring the $i$ prescribed
pairs to occur gives, for $3\le i\le n$,
\[
  \frac{\EE\bigl(\ell_i(G^*)\bigr)}{(d-1)^i}
     =\frac1{2i}\prod_{j=0}^{i-1}
        \frac{dn-dj}{dn-2j-1}.
\]
Write $R_i$ for the product and $N=\lfloor\sqrt n\rfloor$, and assume $n\ge9$.
Each factor is at least $1-j/n$, and every factor with $j\ge1$ is at most one. Hence
\[
  1-\frac{i(i-1)}{2n}\le\prod_{j=0}^{i-1}(1-j/n)
    \le R_i\le\frac{dn}{dn-1}.
\]
For $3\le i\le N$, this gives
$\sum_{i=3}^N |R_i-1|/(2i)\le\sum_{i=3}^N(i-1)/(4n)\le1/8$.
Since $\log\sqrt n\le\sum_{i=1}^N1/i\le1+\log\sqrt n$, it follows that
\[
 \left|\sum_{i=3}^N\frac{R_i}{2i}-\frac14\log n\right|\le\frac78 .
\]
For the tail put $R_{n+1}=0$. The product gives, for $N<i\le n$,
\[
 R_i-R_{i+1}
   =R_i\frac{(d-2)i-1}{dn-2i-1}
   \ge R_i\frac{i-1}{3n},
\]
since $d\ge3$. Summing these differences gives
\[
 \sum_{i=N+1}^n\frac{R_i}{2i}
 \le\frac{3n}{2N(N+1)}R_{N+1}
 \le\frac32\,\frac{N+2}{N+1}\,\frac{dn}{dn-1}<2,
\]
where $n\le N(N+2)$, $N\ge3$ and $dn\ge27$.

By Janson's simplicity formula~\cite[Corollary~1.5]{Janson09}, the
pairing is simple with probability $p_d+o(1)$, where
$p_d=\exp(-(d^2-1)/4)$. Fixing an $i$-cycle reduces the residual
degrees by two at $i$ vertices. For $i\le\sqrt n$ this changes the
degree sums in the formula by $O_d(\sqrt n)$, so the residual pairing
is simple with probability $p_d+o(1)$ uniformly. The probability it
duplicates a prescribed cycle edge is $O_d(i/n)$ by a union bound.
Thus conditioning on simplicity multiplies each expectation through
$\sqrt n$ by $1+o(1)$ uniformly. For all sufficiently large $n$, the
simplicity probability is at least $p_d/2$, so the conditional expected tail
is at most $2/(p_d/2)=4e^{(d^2-1)/4}$. As all simple $d$-regular graphs have the same number
of pairings, their mean weighted cycle sum is
$(\frac14+o(1))\log n$. In particular the leading coefficient in
\eqref{eq:weightedcyclesum} is sharp for each fixed $d\ge3$.

\subsection{The complete graph.}
\label{subsec:complete-graph}
Let $n$ tend to infinity through odd integers. Since
$\ell_i(K_n)=(n)_i/(2i)$ for $3\le i\le n$, where
$(n)_i=n(n-1)\cdots(n-i+1)$, we have
\[
  \ell(K_n)=\frac12\sum_{i=3}^n\frac{1}{i}
       \prod_{j=0}^{i-1}\frac{n-j}{n-2}
     =\frac12\sum_{i=3}^{\infty}\frac{e^{-i^2/(2n)}}{i}+o(1).
\]
Indeed, for $i\le\sqrt n\log n$, the logarithm of the product is
$-i^2/(2n)+O(i/n+i^3/n^2)$ by Taylor expansion. The resulting
error in the sum is $o(1)$. For larger $i$ the inequality
$\log(1+x)\le x$ bounds the product by a constant times
$e^{-i^2/(4n)}$, so both tails are $o(1)$.

Put $a=\sqrt{2n}$. Subtracting the harmonic sum through $\lfloor a\rfloor$
from $\sum_{i\ge1}e^{-(i/a)^2}/i$ gives a Riemann sum tending to
\[
  \int_0^\infty\frac{e^{-x^2}-\mathbf{1}_{\{x\le1\}}}{x}\,dx
  =\frac12\int_0^\infty e^{-t}\log t\,dt=-\frac{\gamma}{2},
\]
by the substitution $t=x^2$ and integration by parts. Restoring the
harmonic sum and omitting $i=1,2$ gives
\[
  \ell(K_n)=\frac14\log(2n)+\frac14\gamma-\frac34+o(1).
\]
For $d=n-1$, the coefficient $\lambda=(n-2)\log(1+1/(n-2))$
is $1+O(1/n)$. Stirling's formula and McKay's enumeration of regular
tournaments~\cite{McK90} give
\[
  \varepsilon(K_n)
      =\left(\sqrt{\frac{\pi}{2}}\,e^{-3/4}+o(1)\right)n\cdot S_{n-1}(n).
\]
It follows that the right side of \eqref{eq:cyclebound} for $K_n$
equals
$\bigl(\sqrt{2/\pi}\,(2e^\gamma)^{1/4}+o(1)\bigr)
 n^{-3/4}\varepsilon(K_n)$.
\section{Proof of Theorem \ref{thm:main2}}\label{sec:multigraph-upper}

\begingroup
\setlength{\parskip}{.42em plus 1pt minus 1pt}
\setlength{\abovedisplayskip}{.45\baselineskip plus 2pt minus 2pt}
\setlength{\belowdisplayskip}{.45\baselineskip plus 2pt minus 2pt}
\setlength{\abovedisplayshortskip}{.2\baselineskip plus 1pt}
\setlength{\belowdisplayshortskip}{.35\baselineskip plus 1pt minus 1pt}
\clubpenalty=10000
\widowpenalty=10000

We proceed by induction on $n$. At each step we remove any components listed below, then reduce triple edges, double edges, and vertices of degree two, in that order. After each reduction we start again with the base cases. Write the upper bound in Theorem~\ref{thm:main2} as
\[
 f(n,m_2,m_3,m_4,k)
 =24^{\frac{n}{5}-\frac{m_2}{5}-\frac{2m_3}{5}-\frac{2m_4}{5}-\frac{k}{10}}
  \cdot10^{\frac{m_2}{3}-\frac{k}{6}}\cdot2^{\frac{k}{2}}\cdot6^{m_3+m_4}.
\]

\textbf{Base cases.} The empty graph has one Eulerian orientation, and $f(0,0,0,0,0)=1$. For each component in the table, $\varepsilon(G)\le f(n,m_2,m_3,m_4,k)$. We may remove any such component, since both the orientation count and $f$ are multiplicative over components.

\begingroup
\centering\small
\setlength{\tabcolsep}{3pt}
\tikzset{v/.append style={inner sep=1.15pt},e/.style={line width=.55pt},de/.style={line width=.55pt}}
\begin{tabular*}{\linewidth}{@{\extracolsep{\fill}}ccccc@{\hspace{1.4em}}ccccc@{}}
\toprule
Case & $G$ & $(n,m_2,m_3,m_4,k)$ & $\varepsilon(G)$ & $f$
& Case & $G$ & $(n,m_2,m_3,m_4,k)$ & $\varepsilon(G)$ & $f$ \\
\midrule
1 &
  \begin{tikzpicture}[scale=1.02,baseline={([yshift=-.55ex]current bounding box.center)}]
    \node[v] (a) at (0, 0.6) {};
    \node[v] (b) at (-0.57, 0.19) {};
    \node[v] (c) at (-0.35, -0.49) {};
    \node[v] (d) at (0.35, -0.49) {};
    \node[v] (e) at (0.57, 0.19) {};

    \draw[e] (a) -- (b);
    \draw[e] (a) -- (c);
    \draw[e] (a) -- (d);
    \draw[e] (a) -- (e);
    \draw[e] (b) -- (c);
    \draw[e] (b) -- (d);
    \draw[e] (b) -- (e);
    \draw[e] (c) -- (d);
    \draw[e] (c) -- (e);
    \draw[e] (d) -- (e);
  \end{tikzpicture} & $(5,0,0,0,0)$ & $24$ & $24$
& 4 &
  \begin{tikzpicture}[scale=1.02,baseline={([yshift=-.55ex]current bounding box.center)}]
    \node[v] (left) at (0, 0) {};
    \node[v] (right) at (1, 0) {};
    \draw[de] (right) to[out=130, in=50] (left);
    \draw[de] (right) to[out=-130, in=-50] (left);
  \end{tikzpicture} & $(2,1,0,0,2)$ & $2$ & $2$ \\[2pt]
2 &
  \begin{tikzpicture}[scale=1.02,baseline={([yshift=-.55ex]current bounding box.center)}]
    \node[v] (top) at (0, 0.35) {};
    \node[v] (left) at (-0.4, -0.3) {};
    \node[v] (right) at (0.4, -0.3) {};
    \draw[de] (top) to[out=-155, in=90] (left);
    \draw[de] (top) to[out=-90, in=30] (left);
    \draw[de] (top) to[out=-25, in=90] (right);
    \draw[de] (top) to[out=-90, in=150] (right);
    \draw[de] (right) to[out=150, in=25] (left);
    \draw[de] (right) to[out=-150, in=-25] (left);
  \end{tikzpicture} & $(3,3,0,0,0)$ & $10$ & $10$
& 5 &
  \begin{tikzpicture}[scale=1.02,baseline={([yshift=-.55ex]current bounding box.center)}]
    \node[v] (a1) at (-0.4,0) {}; 
    \node[v] (b1) at (0.4,0) {}; 
    \node[v] (c1) at (0,0.5) {}; 
    \node[v] (d1) at (0,-0.5) {};
    \draw[de] (a1) to[out=40, in=140] (b1);
    \draw[de] (a1) to[out=-40, in=-140] (b1);
    \draw[e] (a1) -- (c1) -- (b1);
    \draw[e] (a1) -- (d1) -- (b1);
  \end{tikzpicture} & $(4,1,0,0,2)$ & $6$ & $\approx 7.13$ \\[2pt]
3 &
  \begin{tikzpicture}[scale=1.02,baseline={([yshift=-.55ex]current bounding box.center)}]
    \node[v] (left) at (0, 0) {};
    \node[v] (right) at (1, 0) {};
    \draw[de] (right) to[out=130, in=50] (left);
    \draw[de] (right) to[out=160, in=20] (left);
    \draw[de] (right) to[out=-160, in=-20] (left);
    \draw[de] (right) to[out=-130, in=-50] (left);
  \end{tikzpicture} & $(2,0,0,1,0)$ & $6$ & $6$
& 6 &
  \begin{tikzpicture}[scale=1.02,baseline={([yshift=-.55ex]current bounding box.center)}]
    \node[v] (a1) at (-0.4,-0.4) {};
    \node[v] (b1) at (-0.4,0.4) {};
    \node[v] (c1) at (0.4,0.4) {};
    \node[v] (d1) at (0.4,-0.4) {};
    \draw[e] (a1) -- (b1); 
    \draw[e] (a1) -- (c1); 
    \draw[e] (c1) -- (d1); 
    \draw[e] (b1) -- (d1); 
    \draw[de] (b1) to[out=30, in=150] (c1);
    \draw[de] (b1) to[out=-30, in=-150] (c1);
    \draw[de] (a1) to[out=30, in=150] (d1);
    \draw[de] (a1) to[out=-30, in=-150] (d1);
  \end{tikzpicture} & $(4,2,0,0,0)$ & $16$ & $\approx 16.55$ \\[2pt]
\bottomrule
\end{tabular*}\par
\endgroup

For base case 1, fix a vertex of $K_5$. Choose its two out-neighbors and orient the edge between them and the edge between the two in-neighbors. Balance then forces the four edges between the two pairs, giving $\binom42\cdot2^2=24$ orientations. In base case 2 either every double edge is split, giving $2^3$ choices, or all three double edges run around the triangle in one of its two directions. Base cases 3 and 4 have $\binom42$ and $\binom21$ orientations. Suppressing the two degree-two vertices in base case 5 gives base case 3 and preserves the orientation count. For base case 6, exactly two of the four single edges point from one doubled pair to the other. If these two edges form a matching, both doubled pairs split, giving four choices. Each of the other four choices forces one doubled pair and splits the other, giving two choices. This yields $2\cdot4+4\cdot2=16$.

Only base case~3 has an edge of multiplicity four, so no such edge remains after the base cases.

\textbf{Reductions on Triple Edges.}
For a reduction from $G$ to $G'$, let $n',m_2',m_3',m_4',k'$ be the
parameters of $G'$ and write
\[
 f'=f(n',m_2',m_3',m_4',k'),\qquad
 \Delta f=\frac{f(n,m_2,m_3,m_4,k)}{f'},\qquad
 \Delta\varepsilon=\frac{\varepsilon(G)}{\varepsilon(G')}.
\]
By induction $\varepsilon(G')\le f'$, so it suffices to show
$\Delta\varepsilon\le\Delta f$.

Each end of a triple edge has one further incident edge. Delete the triple
edge and identify its ends to give a vertex of degree two. In any Eulerian
orientation of $G'$ one boundary edge enters this vertex and the other
leaves. To extend the orientation to $G$, choose which of the three deleted
edge copies points against the other two. The boundary directions then force
all three directions. The outside ends are distinct in case~1 and coincide
in case~2. If their common neighbor has degree two, its other two incident
edges in the drawing are absent.

\begingroup
\small
\tikzset{every picture/.append style={scale=.84,baseline={([yshift=-.45ex]current bounding box.center)}},
         v/.append style={inner sep=1.15pt},e/.append style={line width=.55pt},
         de/.append style={line width=.55pt}}
\noindent\begin{tabular*}{\linewidth}{@{}l@{\extracolsep{\fill}}cccc@{}}
\toprule
Case & Reduction & $f'$ & $\Delta\varepsilon$ & $\Delta f$ \\
\midrule
1 & \begin{tikzpicture}
    \node[v] (a) at (-0.4,0.2) {}; 
    \node[v] (b) at (0.4,0.2) {};
    \draw[de] (a) to[out=50, in=130] (b);
    \draw[de] (a) to[out=-50, in=-130] (b);
    \draw[e] (a) -- (b);
    \draw[e] (-0.9, -0.3) -- (a);
    \draw[e] (0.9, -0.3) -- (b);
    \node at (2.0,0) {$\longrightarrow$};
    \node[v] (c) at (3.5,0.2) {};
    \draw[e] (3.0,-0.2) -- (c); 
    \draw[e] (4.0,-0.2) -- (c);
    \path (-1.12,0) (4.52,0);
  \end{tikzpicture} & $f(n-1,m_2,m_3-1,m_4,k+1)$ & $3$ & $\approx 4.53$ \\
\addlinespace[2pt]
2 & \begin{tikzpicture}
    \node[v] (a) at (-0.4,0.4) {}; 
    \node[v] (b) at (0.4,0.4) {};
    \node[v] (c) at (0,-0.2) {};
    \draw[de] (a) to[out=50, in=130] (b);
    \draw[de] (a) to[out=-50, in=-130] (b);
    \draw[e] (a) -- (b) -- (c) -- (a);
    \draw[e] (0.4, -0.7) -- (c) -- (-0.4, -0.7);
    \node at (2.0,0) {$\longrightarrow$};
    \node[v] (d) at (3.5,0.4) {};
    \node[v] (e) at (3.5,-0.3) {};
    \draw[de] (d) to[out=-40, in=40] (e);
    \draw[de] (d) to[out=-140, in=140] (e);
    \draw[e] (3.1,-0.7) -- (e); 
    \draw[e] (3.9,-0.7) -- (e);
    \path (-1.12,0) (4.52,0);
  \end{tikzpicture} & $f(n-1,m_2+1,m_3-1,m_4,k+1)$ & $3$ & $\approx 3.97$ \\
\bottomrule
\end{tabular*}\par
\endgroup

\textbf{Reductions on Double Edges.} We check the cases below in order, restarting the induction after each reduction. In the first three cases $u,v$ are the left and right ends of the double edge, and $a,b$ are their top and bottom common neighbors. The outside ends of the edges leaving a drawing may coincide.

\begingroup
\tikzset{double reduction/.style={scale=1.06,
  v/.append style={inner sep=1.1pt},
  e/.append style={line width=.55pt},de/.append style={line width=.55pt}}}

\textbf{1.} Suppose $ab$ is a single edge. Replace the four vertices by a vertex of degree two. Each of $a,b$ has one edge leaving the four vertices. Put $h=1$ if their outside ends coincide and $h=0$ otherwise.
\[
\vcenter{\hbox{%
\begin{tikzpicture}[double reduction]

    \node[v] (a) at (-0.4,0) {}; 
    \node[v] (b) at (0.4,0) {}; 
    \node[v] (c) at (0,0.5) {}; 
    \node[v] (d) at (0,-0.5) {};
    \draw[e] (-0.5,0.8) -- (c);
    \draw[e] (-0.5,-0.8) -- (d);
    \draw[de] (c) -- (d);
    \draw[de] (a) to[out=40, in=140] (b);
    \draw[de] (a) to[out=-40, in=-140] (b);
    \draw[e] (a) -- (c) -- (b);
    \draw[e] (a) -- (d) -- (b);
    \node at (2.0,0) {$\longrightarrow$};
    \node[v] (c2) at (3.5,0) {};
    \draw[e] (3.0,0.3) -- (c2) -- (3.0,-0.3);
  \end{tikzpicture}%
}}\qquad
\begin{aligned}
 f'&=f(n-3,m_2-1+h,m_3,m_4,k+1),\\
 \Delta\varepsilon&=8,\qquad \Delta f\approx\begin{cases}10.95 & (h=0),\\ 9.60 & (h=1).\end{cases}
\end{aligned}
\]
Closing the four vertices by adding a second edge $ab$ gives base case~6, with $16$ orientations. At the replacement vertex one boundary edge enters and one leaves, prescribing the direction of the added edge. Reversal pairs its two directions, so each orientation of $G'$ has eight extensions. When $h=1$, the two boundary edges form a new double edge.

\textbf{2.} Suppose $ab$ is absent and $a,b$ have degrees two and four, respectively, interchanging them if necessary. Replace the four vertices by a vertex of degree two.
\[
\vcenter{\hbox{%
\begin{tikzpicture}[double reduction]

    \node[v] (a) at (-0.4,0) {}; 
    \node[v] (b) at (0.4,0) {}; 
    \node[v] (c) at (0,0.5) {}; 
    \node[v] (d) at (0,-0.5) {};
    \draw[e] (-0.5,-0.8) -- (d) -- (0.5,-0.8);
    \draw[de] (a) to[out=40, in=140] (b);
    \draw[de] (a) to[out=-40, in=-140] (b);
    \draw[e] (a) -- (c) -- (b);
    \draw[e] (a) -- (d) -- (b);
    \node at (2.0,0) {$\longrightarrow$};
    \node[v] (c2) at (3.5,0) {};
    \draw[e] (3.0,-0.4) -- (c2); 
    \draw[e] (4.0,-0.4) -- (c2);
  \end{tikzpicture}%
}}\qquad
\begin{aligned}
 f'&=f(n-3,m_2-1,m_3,m_4,k),\\
 \Delta\varepsilon&=6,\qquad \Delta f\approx7.68.
\end{aligned}
\]
Both boundary edges were incident with $b$, so their multiplicity is preserved even if their outside ends coincide. One enters and one leaves the removed set. Each of $uav$ and $ubv$ is therefore directed consistently along the path. These paths and the two copies of $uv$ give four connections, exactly two directed towards $v$. Thus there are $\binom42=6$ extensions. If $a,b$ both have degree two, the four vertices form base case~5.

\textbf{3.} Suppose $ab$ is absent, $a,b$ have degree four, and they have no common outside neighbor. Replace the four vertices by a vertex of degree four. Every existing double edge from $a$ or $b$ to the outside stays double, and no new multiple edge is formed.
\[
\vcenter{\hbox{%
\begin{tikzpicture}[double reduction]

    \node[v] (a) at (-0.4,0) {}; 
    \node[v] (b) at (0.4,0) {}; 
    \node[v] (c) at (0,0.5) {}; 
    \node[v] (d) at (0,-0.5) {};
    \draw[e] (-0.5,0.8) -- (c) -- (0.5,0.8);
    \draw[e] (-0.5,-0.8) -- (d) -- (0.5,-0.8);
    \draw[de] (a) to[out=40, in=140] (b);
    \draw[de] (a) to[out=-40, in=-140] (b);
    \draw[e] (a) -- (c) -- (b);
    \draw[e] (a) -- (d) -- (b);
    \node at (2.0,0) {$\longrightarrow$};
    \node[v] (c2) at (3.5,0) {};
    \draw[e] (3.0,0.3) -- (c2); 
    \draw[e] (3.0,-0.3) -- (c2); 
    \draw[e] (4.0,0.3) -- (c2); 
    \draw[e] (4.0,-0.3) -- (c2);
  \end{tikzpicture}%
}}\qquad
\begin{aligned}
 f'&=f(n-3,m_2-1,m_3,m_4,k),\\
 \Delta\varepsilon&\le6,\qquad \Delta f\approx7.68.
\end{aligned}
\]
Fix the four boundary directions in an orientation of $G'$. If the two edges at $a$ both enter or both leave, the two at $b$ have the opposite direction. The four internal edges incident with $a,b$ are forced, and the two copies of $uv$ split in two ways. Otherwise each of $a,b$ sends one edge to $u$ or $v$. Sending both edges to the same vertex forces both copies of $uv$, while sending them to different vertices forces the copies to split. This gives $2\cdot1+2\cdot2=6$ extensions.

\textbf{3(a).} Suppose instead that $a,b$ have a common outside neighbor. Delete $u$ and its incident edges. To form $G_0$ add $ab$, and to form $G_1$ add one copy each of $va$ and $vb$. Both reductions retain every edge outside the displayed configuration. The drawing shows two common outside neighbors $x,y$.
\[
\vcenter{\hbox{%
\begin{tikzpicture}[scale=1.30,v/.append style={inner sep=1.35pt},baseline=-0.5ex,every label/.style={font=\scriptsize,inner sep=1.25pt}]
  \begin{scope}
    \node[v,label=above:$a$] (a) at (0,.68) {};
    \node[v,label=below:$b$] (b) at (0,-.68) {};
    \node[v,label=left:$x$] (x) at (-1.02,0) {};
    \node[v,label=right:$y$] (y) at (1.02,0) {};
    \node[v,label=below left:$u$] (u) at (-.33,0) {};
    \node[v,label=below right:$v$] (v) at (.33,0) {};
    \foreach \w in {x,u,v,y} \draw[e] (a)--(\w)--(b);
    \draw[de] (u) to[bend left=24] (v);
    \draw[de] (u) to[bend right=24] (v);
    \node at (0,-1.23) {$G$};
  \end{scope}
  \node at (1.64,0) {$\longrightarrow$};
  \begin{scope}[xshift=3.3cm]
    \node[v,label=above:$a$] (a0) at (0,.68) {};
    \node[v,label=below:$b$] (b0) at (0,-.68) {};
    \node[v,label=left:$x$] (x0) at (-1.02,0) {};
    \node[v,label=right:$y$] (y0) at (1.02,0) {};
    \node[v,label=right:$v$] (v0) at (.4,0) {};
    \foreach \w in {x0,v0,y0} \draw[e] (a0)--(\w)--(b0);
    \draw[e] (a0)--(b0);
    \node at (0,-1.23) {$G_0$};
  \end{scope}
  \node at (4.88,0) {\small and};
  \begin{scope}[xshift=6.5cm]
    \node[v,label=above:$a$] (a1) at (0,.68) {};
    \node[v,label=below:$b$] (b1) at (0,-.68) {};
    \node[v,label=left:$x$] (x1) at (-1.02,0) {};
    \node[v,label=right:$y$] (y1) at (1.02,0) {};
    \node[v,label=right:$v$] (v1) at (0,0) {};
    \foreach \w in {x1,y1} \draw[e] (a1)--(\w)--(b1);
    \draw[de] (a1) to[bend left=24] (v1) to[bend left=24] (b1);
    \draw[de] (a1) to[bend right=24] (v1) to[bend right=24] (b1);
    \node at (0,-1.23) {$G_1$};
  \end{scope}
\end{tikzpicture}%
}}
\]
Each Eulerian orientation pairs the two incoming edges at $u$ with the two outgoing edges in two ways. Pairing the copies of $uv$ together gives $G_0$ with a loop at $v$, which has two possible directions. Each of the other two pairings gives $G_1$. Hence $2\varepsilon(G)=2\varepsilon(G_0)+2\varepsilon(G_1)$. In $G_0$, $v$ has degree two and the new edge $ab$ is single. In $G_1$, $va,vb$ are doubled and $uv$ is removed. For any outside ends, induction gives
\begin{align*}
 \varepsilon(G)
 &\le f(n-1,m_2-1,m_3,m_4,k+1)+f(n-1,m_2+1,m_3,m_4,k)\\
 &=\left(\frac{\sqrt2}{24^{1/10}\sqrt{10}}+
              \frac{10^{1/3}}{24^{2/5}}\right)f(n,m_2,m_3,m_4,k)
 <0.931\,f(n,m_2,m_3,m_4,k).
\end{align*}

\textbf{4.} Suppose $u,v$ have exactly one common neighbor $w$, both $uw$ and $vw$ are single, and $w$ has degree four. Contract the triangle to a vertex of degree four, keeping the four boundary edges. For distinct outside ends, as drawn, we have
\[
\vcenter{\hbox{%
\begin{tikzpicture}[double reduction]

    \node[v] (a) at (-0.4,-0.2) {}; 
    \node[v] (b) at (0.4,-0.2) {}; 
    \node[v] (c) at (0,0.3) {};
    \draw[e] (a) -- (-1,-0.5); 
    \draw[e] (b) -- (1,-0.5);
    \draw[e] (-0.5,0.6) -- (c) -- (0.5,0.6);
    \draw[de] (a) to[out=40, in=140] (b);
    \draw[de] (a) to[out=-40, in=-140] (b);
    \draw[e] (a) -- (c) -- (b);
    \node[below=1.5pt,font=\scriptsize] at (a) {$u$};
    \node[below=1.5pt,font=\scriptsize] at (b) {$v$};
    \node[above=1.5pt,font=\scriptsize] at (c) {$w$};
    \node at (2.0,0) {$\longrightarrow$};
    \node[v] (c2) at (3.5,0) {};
    \draw[e] (3.0,0.3) -- (c2); 
    \draw[e] (3.0,-0.3) -- (c2); 
    \draw[e] (4.0,0.3) -- (c2); 
    \draw[e] (4.0,-0.3) -- (c2);
  \end{tikzpicture}%
}}\qquad
\begin{aligned}
 f'&=f(n-2,m_2-1,m_3,m_4,k),\\
 \Delta\varepsilon&\le3,\qquad \Delta f\approx4.07.
\end{aligned}
\]
Fix an orientation of $G'$. If the boundary edges at $u,v$ both enter the triangle or both leave, $uw,vw$ are forced and $uv$ splits in two ways. If one enters and the other leaves, precisely one of $wu,wv$ points away from $w$. One choice forces both copies of $uv$ and the other requires them to split, giving three extensions.

Write $x,y$ for the outside ends at $u,v$, and $z,z'$ for those at $w$. Since $w$ is the only common neighbor, $x\ne y$. If no three ends coincide, let $s$ count the new double edges formed by an edge from $\{ux,vy\}$ and one from $\{wz,wz'\}$. An existing double edge with $z=z'$ stays double. Thus $0\le s\le2$ and
\[
 f'=f(n-2,m_2-1+s,m_3,m_4,k),\qquad
 \Delta f\ge \frac{24^{3/5}}{10^{1/3}}>3.
\]
If three ends coincide, say $x=z=z'$, the double edges $uv,wx$ are replaced by a triple edge. Fix the directions of the other edges in $G'$. Exactly one of the three copies at $x$ agrees in direction with $vy$ at the contracted vertex. In one of the three choices the boundary directions at $u,v$ agree, and in two they differ. The three orientations therefore have $2+3+3=8$ extensions in total. Hence
\[
 \Delta\varepsilon=\frac83,\qquad f'=f(n-2,m_2-2,m_3+1,m_4,k),\qquad
 \Delta f=\frac{24^{2/5}10^{2/3}}6>\frac83.
\]

\textbf{5.} Suppose $u,v$ have exactly one common neighbor $w$ and $w$ has degree two. Contract the triangle to a vertex of degree two. Its two outside ends are distinct, and the two boundary directions are opposite. Case~4 then gives three extensions.
\[
\vcenter{\hbox{%
\begin{tikzpicture}[double reduction]

    \node[v] (a) at (-0.4,0) {}; 
    \node[v] (b) at (0.4,0) {}; 
    \node[v] (c) at (0,0.5) {};
    \draw[e] (a) -- (-1,-0.3); 
    \draw[e] (b) -- (1,-0.3);
    \draw[de] (a) to[out=40, in=140] (b);
    \draw[de] (a) to[out=-40, in=-140] (b);
    \draw[e] (a) -- (c) -- (b);
    \node[below=1.5pt,font=\scriptsize] at (a) {$u$};
    \node[below=1.5pt,font=\scriptsize] at (b) {$v$};
    \node[above=1.5pt,font=\scriptsize] at (c) {$w$};
    \node at (2.0,0) {$\longrightarrow$};
    \node[v] (d) at (3.5,0) {};
    \draw[e] (3.0,-0.4) -- (d); 
    \draw[e] (4.0,-0.4) -- (d);
  \end{tikzpicture}%
}}\qquad
\begin{aligned}
 f'&=f(n-2,m_2-1,m_3,m_4,k),\\
 \Delta\varepsilon&=3,\qquad \Delta f\approx4.07.
\end{aligned}
\]

\textbf{5(a).} Suppose $u,v$ have exactly one common neighbor, joined to one endpoint by a second double edge. Three doubled sides give base case~2. Otherwise two sides are doubled and one is single. Contract this triangle to a vertex of degree two and let $h=1$ if the outside ends coincide and $h=0$ otherwise.
\[
\vcenter{\hbox{%
\begin{tikzpicture}[double reduction]

    \node[v] (a) at (-0.4,-0.22) {};
    \node[v] (b) at (0.4,-0.22) {};
    \node[v] (c) at (0,0.38) {};
    \draw[e] (-.92,-.48) -- (a) -- (b) -- (.92,-.48);
    \draw[de] (a) to[bend left=17] (c);
    \draw[de] (a) to[bend right=17] (c);
    \draw[de] (b) to[bend left=17] (c);
    \draw[de] (b) to[bend right=17] (c);
    \node at (2.0,0) {$\longrightarrow$};
    \node[v] (d) at (3.5,0.08) {};
    \draw[e] (3.0,-.32) -- (d) -- (4.0,-.32);
  \end{tikzpicture}%
}}\qquad
\begin{aligned}
 f'&=f(n-2,m_2-2+h,m_3,m_4,k+1),\\
 \Delta\varepsilon&=5,\qquad \Delta f\ge24^{3/10}\sqrt5>5.
\end{aligned}
\]
The boundary edges meet the endpoints of the single side. Adding a distinguished second copy of this side closes the triangle to base case~2. Reversal pairs its ten orientations according to the direction of the distinguished edge, so each orientation of $G'$ has five extensions.

\textbf{6.} If one endpoint of the double edge has degree two, the other has degree four, by base case~4. Delete the degree-two endpoint and both copies. In each extension one copy is directed each way. A double edge among the other two incidences stays double.
\[
\vcenter{\hbox{%
\begin{tikzpicture}[double reduction]

    \node[v] (a) at (-0.1,0) {}; 
    \node[v] (b) at (0.7,0) {};
    \draw[de] (a) to[out=45, in=135] (b);
    \draw[de] (a) to[out=-45, in=-135] (b);
    \draw[e] (a) -- (-0.6,0.3); 
    \draw[e] (a) -- (-0.6,-0.3);
    \node at (2.0,0) {$\longrightarrow$};
    \node[v] (c) at (3.5,0) {};
    \draw[e] (3.0,0.3) -- (c) -- (3.0,-0.3);
  \end{tikzpicture}%
}}\qquad
\begin{aligned}
 f'&=f(n-1,m_2-1,m_3,m_4,k),\\
 \Delta\varepsilon&=2,\qquad \Delta f=10^{1/3}>2.
\end{aligned}
\]

\textbf{7.} Otherwise both endpoints have degree four and no common neighbor. Delete both copies and identify the endpoints. Their outside neighbors are disjoint, so no new multiple edge is formed. In an orientation of $G'$, the two other incidences at either endpoint determine how many deleted copies must point towards it, giving at most two extensions.
\[
\vcenter{\hbox{%
\begin{tikzpicture}[double reduction]

    \node[v] (a) at (-0.4,0) {}; 
    \node[v] (b) at (0.4,0) {};
    \draw[de] (a) to[out=50, in=130] (b);
    \draw[de] (a) to[out=-50, in=-130] (b);
    \draw[e] (-0.9,0.3) -- (a) -- (-0.9,-0.3);
    \draw[e] (0.9,0.3) -- (b) -- (0.9,-0.3);
    \node at (2.0,0) {$\longrightarrow$};
    \node[v] (c) at (3.5,0) {};
    \draw[e] (3.0,0.3) -- (c); 
    \draw[e] (3.0,-0.3) -- (c); 
    \draw[e] (4.0,0.3) -- (c); 
    \draw[e] (4.0,-0.3) -- (c);
  \end{tikzpicture}%
}}\qquad
\begin{aligned}
 f'&=f(n-1,m_2-1,m_3,m_4,k),\\
 \Delta\varepsilon&\le2,\qquad \Delta f=10^{1/3}>2.
\end{aligned}
\]
\endgroup

These cases cover every double edge. Case~6 handles a degree-two endpoint. Two degree-four endpoints have at most two common neighbors. If they have two, all four edges to those neighbors are single. An edge of multiplicity two between the common neighbors gives base case~6, a single edge gives case~1, and no edge gives base case~5 or case~2, 3 or~3(a), according to their degrees and outside neighbors. For one common neighbor its degree and the two incident multiplicities give case~4, 5 or~5(a). With no common neighbor, case~7 applies.

\textbf{Reductions on vertices of degree two.}
After all multiple edges have been resolved, let $v$ be a degree-two vertex
with distinct neighbors $a,b$. Replace the path $avb$ by an edge $ab$.
In an Eulerian orientation $avb$ is directed from one end to the other,
so this gives a bijection of orientations and preserves the degrees of
$a$ and $b$. In case~1, $a,b$ are nonadjacent, regardless of their degrees.
In cases~2 and~3 they are adjacent and the new edge creates a double edge.
Case~2 shows a triangle component. Case~3 also applies if just one of
$a,b$ has degree two, with the two outside edges at that vertex omitted.

\begingroup
\small
\tikzset{every picture/.append style={scale=.84,baseline={([yshift=-.45ex]current bounding box.center)}},
         v/.append style={inner sep=1.15pt},e/.append style={line width=.55pt},
         de/.append style={line width=.55pt}}
\noindent\begin{tabular*}{\linewidth}{@{}l@{\extracolsep{\fill}}cccc@{}}
\toprule
Case & Reduction & $f'$ & $\Delta\varepsilon$ & $\Delta f$ \\
\midrule
1 & \begin{tikzpicture}
    \node[v] (left) at (-0.5, 0) {};
    \node[v] (top) at (0, 0.35) {};
    \node[v] (right) at (0.5, 0) {};
    \draw[e] (left) -- (top) -- (right);
    \draw[e] (left) -- (-1.0, 0.4);
    \draw[e] (left) -- (-1.1, 0.0);
    \draw[e] (left) -- (-0.9, -0.35);
    \draw[e] (right) -- (1.0, 0.4);
    \draw[e] (right) -- (1.1, 0.0);
    \draw[e] (right) -- (0.9, -0.35);
    \node at (2.0,0) {$\longrightarrow$};
    \node[v] (r1) at (3.1, 0) {};
    \node[v] (r2) at (3.9, 0) {};
    \draw[e] (r1) -- (r2);
    \draw[e] (r1) -- (2.6, 0.4);
    \draw[e] (r1) -- (2.5, 0.0);
    \draw[e] (r1) -- (2.7, -0.35);
    \draw[e] (r2) -- (4.4, 0.4);
    \draw[e] (r2) -- (4.5, 0.0);
    \draw[e] (r2) -- (4.3, -0.35);
    \path (-1.12,0) (4.52,0);
  \end{tikzpicture} & $f(n-1,m_2,m_3,m_4,k-1)$ & $1$ & $\approx 1.32$ \\
\addlinespace[2pt]
2 & \begin{tikzpicture}
    \node[v] (top) at (0, 0.35) {};
    \node[v] (left) at (-0.4, -0.3) {};
    \node[v] (right) at (0.4, -0.3) {};
    \draw[e] (left) -- (top) -- (right) -- (left);
    \node at (2.0,0) {$\longrightarrow$};
    \node[v] (r1) at (3.1, 0) {};
    \node[v] (r2) at (3.9, 0) {};
    \draw[de] (r1) to[out=45, in=135] (r2);
    \draw[de] (r1) to[out=-45, in=-135] (r2);
    \path (-1.12,0) (4.52,0);
  \end{tikzpicture} & $f(n-1,m_2+1,m_3,m_4,k-1)$ & $1$ & $\approx 1.16$ \\
\addlinespace[2pt]
3 & \begin{tikzpicture}
    \node[v] (top) at (0, 0.45) {};
    \node[v] (left) at (-0.5, 0) {};
    \node[v] (right) at (0.5, 0) {};
    \draw[e] (left) -- (top) -- (right) -- (left);
    \draw[e] (left) -- (-1.0, 0.35);
    \draw[e] (left) -- (-1.0, -0.35);
    \draw[e] (right) -- (1.0, 0.35);
    \draw[e] (right) -- (1.0, -0.35);
    \node at (2.0,0) {$\longrightarrow$};
    \node[v] (r1) at (3.1, 0) {};
    \node[v] (r2) at (3.9, 0) {};
    \draw[de] (r1) to[out=50, in=130] (r2);
    \draw[de] (r1) to[out=-50, in=-130] (r2);
    \draw[e] (r1) -- (2.6, 0.35);
    \draw[e] (r1) -- (2.6, -0.35);
    \draw[e] (r2) -- (4.4, 0.35);
    \draw[e] (r2) -- (4.4, -0.35);
    \path (-1.12,0) (4.52,0);
  \end{tikzpicture} & $f(n-1,m_2+1,m_3,m_4,k-1)$ & $1$ & $\approx 1.16$ \\
\bottomrule
\end{tabular*}\par
\endgroup
\par\endgroup

\textbf{Inductive Step for Simple Graphs.} We consider the local reduction at a vertex $v \in V(G)$ of degree $d(v) = 4$ in a graph $G$ with no vertices of degree $2$ and no double edges. Let $t$ denote the number of triangles in $G$ containing $v$. Since there are $\binom{4}{2} = 6$ potential pairs of neighbors in $N(v)$, $0 \le t \le 6$.

To perform the reduction, we delete $v$ and add two edges in $N(v)$ to pair its four neighbors into a pairing. There are three pairings of $N(v)$. Every balanced orientation at $v$ is compatible with exactly two of them, the two bijections of its incoming edges with its outgoing edges. For a fixed pairing, an orientation after smoothing extends uniquely over $v$. Denote the graph from a pairing $P$ by $G_P$. Then $2\varepsilon(G)=\sum_P\varepsilon(G_P)$, and by induction
\[
  \varepsilon(G)\le \frac12\sum_P f(n-1,m_2'(P),m_3,m_4,k),
\]
where $m_2'(P)$ equals $m_2$, $m_2+1$, or $m_2+2$ according to the number of pairs in $P$ that are already edges of $G$.

For $0\le t\le5$, some two neighbors of $v$ are nonadjacent. The pairing joining those two creates at most one double edge and the other two pairings each create at most two. Since $f$ increases with $m_2$, the right side is greatest when $t=5$. In this case one pairing yields one double edge and the other two pairings yield two double edges. Thus
\begin{align*}
 \varepsilon(G)
 &\le \tfrac12\bigl(f(n-1,m_2+1,m_3,m_4,k)
                    +2f(n-1,m_2+2,m_3,m_4,k)\bigr)\\
 &=\tfrac12\bigl(24^{-2/5}10^{1/3}+2\cdot24^{-3/5}10^{2/3}\bigr)
                   f(n,m_2,m_3,m_4,k)\\
 &<0.992\,f(n,m_2,m_3,m_4,k).
\end{align*}

When $t=6$, all six pairs of neighbors in $N(v)$ are adjacent, so $\{v\}\cup N(v)$ induces an isolated $K_5$. We remove the entire component as in base case 1. Every reduction apart from removing an equality component gives a strict inequality. The first four base cases therefore give exactly the components of graphs attaining equality. In particular, a simple $4$-regular graph attains equality precisely when each component is $K_5$. This completes the proof of Theorem~\ref{thm:main2}. \qed
\section{Random regular graphs}
\label{sec:large-girth}

Creed~\cite[Section~3.4]{Creed10} computes the orientation-weighted short-cycle moments for fixed lengths. 
For a growing cutoff, the proof uses the matching estimate of Lu and Sz\'ekely
\cite[Theorem~5]{LuSzekelyMatching}.

\begin{theorem}\label{thm:eo-conditioned-girth}
Fix an even integer $d\ge4$ and let $G_{n,d}$
be a uniformly random simple $d$-regular graph on $\{1,\ldots,n\}$.
Let $H_r=\sum_{j=1}^r1/j$, and recall $S_d(n)$ from \eqref{eq:Schrijverboundregular}.
If $r=r(n)\ge2$ is an integer and $\Delta_n=r^6(d-1)^{2r}/n=o(1)$, then
\[
  \EE\!\left(\varepsilon(G_{n,d})\,\middle|\,
    \operatorname{girth}(G_{n,d})>r\right)
  = \sqrt{\frac{\pi dn}{2}}\,
    \exp\!\left(-\frac{H_r}{2}+O_d(\Delta_n)\right)\cdot S_d(n).
\]
In particular, for
\[
  r_n=\left\lfloor
       \frac{\log n-8\log(\log n)}{2\log(d-1)}
      \right\rfloor,
\]
the conditional mean is
\[
  \bigl(\sqrt{\pi d\,e^{-\gamma}}+o(1)\bigr)
  \sqrt{\frac{n}{\log_{d-1}n}}\cdot S_d(n),
\]
where $\gamma$ is Euler's constant. For every sufficiently large $n$,
some simple $d$-regular graph of girth greater than $r_n$ has at most
this many Eulerian orientations.
\end{theorem}

\begin{proof}
Put $k=d/2$ and $m=kn$, and give each vertex $d$ labelled half-edges.
We consider two probability distributions, $\PP_1$ and $\PP_2$, on
the set of perfect matchings of these $2m$ half-edges. Each matching
gives a pseudograph by joining the vertices of each matched pair.
Write $\EE_1$ and $\EE_2$ for expectation under the respective distributions.

Under $\PP_1$, choose a perfect matching uniformly at random.
Under $\PP_2$, independently at each vertex choose $k$ of its $d$
half-edges uniformly to be outgoing, and mark the other $k$ as incoming.
Choose a uniform bijection from the $m$ outgoing half-edges to the $m$
incoming half-edges, pair each half-edge with its image, and then forget
the marks. Permuting the half-edge labels at each vertex takes any choice
of marks to any other. Thus, when counting cycles or computing girth
probabilities under $\PP_2$, we may fix the marks in advance.

Let $Y$ count Eulerian orientations of the pairing, counting the two
directions of a loop separately. Counting a pairing together with an
Eulerian orientation gives, for every event $A$ of the underlying graph,
\begin{align*}
  \EE_1(Y\mathbf{1}_A)
    &=\EE_1(Y)\,\PP_2(A),\\
  \EE_1(Y)
    &=\frac{\binom{d}{k}^{n}m!}{(2m-1)!!}
      =\frac{4^m}{\binom{2m}{m}}\cdot S_d(n)
      =\sqrt{\pi m}\,(1+O_d(n^{-1}))\cdot S_d(n).
\end{align*}
Let $X_j=\ell_j(G^*)$, where $G^*$ is the underlying pseudograph. Thus $X_1$ counts loops and $X_2$ counts unordered pairs of parallel edges. Set $A_r=\{X_1=\cdots=X_r=0\}$. As $r\ge2$, every pairing
in $A_r$ gives a simple graph, and every such graph of girth greater
than $r$ has $(d!)^n$ pairings. Thus the desired mean is
$\EE_1(Y)\,\PP_2(A_r)/\PP_1(A_r)$.

Using the falling factorial $(a)_j=a(a-1)\cdots(a-j+1)$, counting cycles with a starting
point and a direction, then dividing by $2j$, gives
\begin{align*}
  \EE_1(X_j)
    &=\frac{(n)_j}{2j}
      \frac{\bigl(d(d-1)\bigr)^j}{\prod_{u=0}^{j-1}(dn-2u-1)}
      =\frac{(d-1)^j}{2j}\bigl(1+O_d(j^2/n)\bigr),\\
  \EE_2(X_j)
    &=\frac{(n)_j}{2j}\frac{\bigl(k(d-1)\bigr)^j+k^j}{(kn)_j}
      =\frac{(d-1)^j+1}{2j}\bigl(1+O_d(j^2/n)\bigr).
\end{align*}
For the second count, record whether each edge points along or against
the traversal. At a vertex, consecutive signs that agree give $k^2$
choices of half-edges, and signs that differ give $k(k-1)$. Summing
over cyclic sign sequences gives the trace of the $j$th power of the
matrix with diagonal entries $k^2$ and off-diagonal entries $k(k-1)$.
Its eigenvalues are $k(d-1)$ and $k$.

We keep track of absolute constants for Corollary~\ref{cor:finite-girth}.
Assume $16r^2\le n$, as holds eventually under the hypothesis of the theorem.
For each $j\le r$, the ratios of the two exact means above to $(d-1)^j/(2j)$
and $((d-1)^j+1)/(2j)$, respectively, have logarithms of absolute value at most
$2j^2/n$. Indeed, the numerator factors in both ratios are $1-u/n$,
and the denominator factors are $1-(2u+1)/(dn)$ and $1-u/(kn)$,
respectively, for $0\le u<j$. Use $|\log(1-x)|\le2x$ for $0\le x\le1/2$
and $d\ge4$. Thus both relative errors are at most $4j^2/n$.
For $b\in\{1,2\}$, write $\mu_b=\sum_{j=1}^r\EE_b(X_j)$. We obtain
\[
 \mu_b\le\sum_{j=1}^r j\,\EE_b(X_j)\le2(d-1)^r,\qquad
 \left|\mu_2-\mu_1-\frac{H_r}{2}\right|\le\frac{8r (d-1)^r}{n}.
\]
Here we used $\sum_{j=1}^r((d-1)^j+1)\le2(d-1)^r$ for $d\ge4$.

We apply the cited matching theorem to the partial matchings specifying
individual cycles of length at most $r$, in $K_{dn}$ for
$\PP_1$ and $K_{m,m}$ for $\PP_2$.
None properly contains another. Each has probability at most $1/n$, since its occurrence requires a specified pair.
By symmetry, the sum of occurrence probabilities over cycles using a
specified half-edge is
\[
  \frac{2}{dn}\sum_{j=1}^r j\,\EE_b(X_j)\le \frac{(d-1)^r}{n}.
\]
This also bounds the sum for a specified pair. Cycles conflicting with
a given cycle $F$ use at least one of its at most $2r$ half-edges,
so their total probability is at most $2r (d-1)^r/n$.

For the remaining hypothesis, condition on the pairs of $F$, and put
$f=|F|$. A distinct $j$-cycle sharing a pair with $F$ can be
traversed starting along one of the $2f$ oriented edges of $F$.
Order the half-edges at each vertex. Each further step selects one of
the $d-1$ half-edges other than the arrival half-edge, so there are at
most $2f (d-1)^{j-1}$ instructions for such a traversal. Follow each
instruction, revealing partners as needed and stopping if an edge
repeats. A successful traversal must use an edge outside $F$, and
its last such edge ends in $V(F)$. When that edge is traversed its
partner is freshly exposed. At any fresh exposure there are at least
$dn/4$ eligible partners, since conditioning on $F$ and following the
traversal exposes at most $2r$ pairs in total. At most $df$ eligible
partners lie in $V(F)$. A union bound over the at most $j$ exposures gives
a success probability at most $4jf/n$. The conditional expected number
of other cycles sharing a pair with $F$ is therefore at most
\[
 \sum_{j=1}^r 2f (d-1)^{j-1}\frac{4jf}{n}
 \le \frac{8r^3(d-1)^r}{(d-2)n}\le\frac{4r^3(d-1)^r}{n}.
\]
This bounds the sum of the probabilities of the residual matchings
$F'\setminus F$ for compatible cycles $F'\ne F$ sharing a pair
with $F$, after the $f$ pairs of $F$ are removed. The matching
theorem evaluates these probabilities after removing $r$ pairs.
Each residual matching specifies at most $r$ pairs. In each factor of its
occurrence probability, the number of eligible partners decreases by at most
$2r$ and remains at least $dn/4$. Thus each residual probability increases
by a factor at most $(1+8r/(dn))^r\le \exp(8r^2/(dn))<2$.

All four probability bounds required by the matching theorem are therefore
strictly smaller than $\eta=16r^3(d-1)^r/n$. If $\eta\le1/64$, the finite
inequalities in the proof of Lu and Sz\'ekely~\cite[Theorem~5]{LuSzekelyMatching}
give
\[
 \left|\log\PP_b(A_r)+\mu_b\right|\le20r\eta\mu_b.
\]
To see the constant, let $y(s)\in[1,2]$ solve $\log y(s)=s y(s)$.
Their lower bound for $\log\PP_b(A_r)$ is
$-\mu_b(y(2\eta)+\eta y(2\eta)^2)\ge-\mu_b(1+9\eta)$, since
$y(2\eta)\le e^{4\eta}\le1+5\eta$.
Their upper bound is $-\mu_b e^{-a}$, where
\[
 a=2\eta y(4\eta)+4\eta^2 y(4\eta)^2+2r\log y(4\eta)
   \le(5+16r)\eta\le20r\eta.
\]
Now use $e^{-a}\ge1-a$. Combining the estimates for $b=1,2$ with the
bounds on the cycle means gives
\begin{equation}\label{eq:finite-girth-ratio}
 \left|\log\frac{\PP_2(A_r)}{\PP_1(A_r)}+\frac{H_r}{2}\right|
 \le20r\eta(\mu_1+\mu_2)+\frac{8r (d-1)^r}{n}
 \le\frac{2048r^4(d-1)^{2r}}{n}.
\end{equation}
Under the hypothesis $\Delta_n=o(1)$, both $16r^2\le n$ and
$\eta\le1/64$ hold for all sufficiently large $n$.
Equation~\eqref{eq:finite-girth-ratio} and the formula for $\EE_1(Y)$
give the first assertion.
For $r=r_n$, we have $\Delta_n=O_d((\log n)^{-2})$ and
$r_n\sim\tfrac12\log_{d-1}n$. Substituting
$H_{r_n}=\log r_n+\gamma+o(1)$ yields the stated constant.
The existence assertion follows by averaging.
\end{proof}

\begin{proof}[Proof of Corollary~\ref{cor:finite-girth}]
Write $t=\sqrt{\log_{d-1} n}$ and
$r=\lfloor(t^2-3t)/2\rfloor$. Then $t\ge8$, $r\ge20$ and
$(d-1)^{2r}/n\le (d-1)^{-3t}$. The hypotheses of \eqref{eq:finite-girth-ratio}
hold, since $16r^2\le4t^4\le3^{t^2}\le n$ and
\[
  \eta=\frac{16r^3(d-1)^r}{n}\le\frac{2t^6}{3^{t^2/2}}<\frac1{64}
  \qquad(t\ge8).
\]
The Wallis bound $4^m/\binom{2m}{m}\le\sqrt{\pi(m+1/2)}$
and $H_r\ge\log r+\gamma$ give, with $m=dn/2$,
\[
 \EE_1(Y\mid A_r)
 \le\frac{\sqrt{\pi d\,e^{-\gamma}n}}{\sqrt{2r}}\,e^{E}\cdot S_d(n),
 \qquad
 E=\frac{2048r^4(d-1)^{2r}}{n}+\frac{1}{2dn}.
\]
Since $r\le t^2/2$ and $d\ge4$,
\[
 E\le 128t^8\,3^{-3t}+\frac{1}{8\cdot3^{t^2}}\le\frac{1}{8t}.
\]
For the last inequality, multiply by $t$. Both terms decrease for
$t\ge8$, and their sum at $t=8$ is less than $1/8$.
Since $t^2-3t-2\le2r\le t^2$, it follows that
\[
 2r e^{-2E}\ge 2r(1-2E)\ge t^2-3t-2-\frac{t}{4}
        \ge (t-2)^2.
\]
The conditional mean is therefore at most
$(\pi d\,e^{-\gamma} n)^{1/2}/(t - 2) \cdot S_d(n)$.
Some simple graph in $A_r$ has at most this many orientations.
Finally, $t/(t-2)\le4/3$ for $t\ge8$, and
$(4/3)(\pi e^{-\gamma})^{1/2} < 2$.
\end{proof}

\section{Concluding remarks}

$\bullet$
For fixed even $d\ge4$, the upper bound in Theorem~\ref{thm:eo-conditioned-girth}
is $\sqrt{2d}+o(1)$ times the conjectured lower bound in Conjecture \ref{conj:lowerbound}. A first step would
be to obtain the factor $(n/\log_{d-1}n)^{1/2}$ with any positive constant
depending only on $d$. The examples in
Section~\ref{subsec:cycle-sharpness} have
$\ell(G)=(1/4+o(1))\log n$, so on these graphs the correction supplied by
Theorem~\ref{thm:main} is only $n^{\lambda/4+o(1)}$.

$\bullet$
Theorem~\ref{thm:main2} gives $\varepsilon(G)\le 24^{n/5}$ for $d=4$,
with equality precisely for disjoint unions of $K_5$, and more generally solves Schrijver's Conjecture -- Conjecture \ref{conj:schrijver-degrees} -- 
for graphs of maximum degree four.  For $d=6$ the
question is whether $\varepsilon(G)\le\varepsilon(K_7)^{n/7}$ for every
simple $6$-regular graph. Disjoint copies of $K_7$ attain this bound when
$7$ divides $n$.

$\bullet$ It would be interesting to know the correct analog of the weighted cycle sum for a lower bound on the number of Eulerian orientations in pseudographs which are not necessarily regular.

$\bullet$
Together with Theorem~\ref{thm:main}, Conjecture~\ref{conj:cycle-upper} would bound
$\log\bigl(\varepsilon(G)/S_d(n)\bigr)$ between $\lambda\ell(G)$ and
$c_d\ell(G)$. The example of two vertices joined by $d$ parallel edges
in Section~1.2 forces
\[
 c_d\ge\frac{2(d-1)}{d}\log\!\Bigl(\frac{2^d}{\binom{d}{d/2}}\Bigr)
     =\log(\pi d/2)+o(1)
 \qquad (d\to\infty).
\]

\section*{AI Declaration}
Generative AI assisted with literature searches, production of figures, and checking and revising proofs. The authors take full responsibility for the material appearing in this paper.

\end{document}